\documentclass{gtart_a}
\pdfoutput=1
\usepackage{colortbl}
\usepackage{pinlabel}

\title{Automorphic forms and rational homology 3--spheres}

\author{Frank Calegari}
\givenname{Frank}
\surname{Calegari}
\address{Department of Mathematics\\
Harvard University\\\newline
One Oxford Street\\
Cambridge MA 02138\\
USA}
\email{fcale@math.harvard.edu}
\urladdr{http://www.math.harvard.edu/~fcale/}

\author{Nathan M Dunfield}
\givenname{Nathan M}
\surname{Dunfield}
\address{Mathematics 253-37\\
Caltech\\\newline
Pasadena CA 91125\\
USA}
\email{dunfield@caltech.edu}
\urladdr{http://www.its.caltech.edu/~dunfield/}

\volumenumber{10}
\issuenumber{}
\publicationyear{2006}
\papernumber{8}
\lognumber{0643}
\startpage{295}
\endpage{329}

\doi{10.2140/gt.2006.10.295}
\MR{}
\Zbl{}


\received{18 August 2005}
\revised{28 February 2006}
\accepted{2 January 2006}
\published{2 April 2006}
\publishedonline{2 April 2006}
\proposed{Walter Neumann}
\seconded{David Gabai, Tomasz Mrowka}
\corresponding{}
\editor{}
\version{}

\AtBeginDocument{\let\bar\wbar}

\swapnumbers
\def\setobjecttype#1{}

\newcommand{\osl}[1]{{\mkern 2mu\overline{\mkern-2mu #1\mkern-1mu}\mkern1mu}}
\newcommand{\ol}[1]{{\mkern 2mu\overline{\mkern-2mu #1\mkern-2mu}\mkern2mu}}

\def\iijj{i\kern-0.1em{j}}
\def\GSp{\mathrm{GSp}}
\def\PGL{\mathrm{PGL}}
\def\PSL{\mathrm{PSL}}
\def\SL{\mathrm{SL}}
\def\GL{\mathrm{GL}}
\def\Ext{\mathrm{Ext}}
\def\Gal{\mathrm{Gal}}
\def\CL{\mathrm{CL}}
\def\twpi{\widetilde{\pi}}
\def\pibar{\osl{\pi}}
\def\rhobar{\osl{\rho}}
\def\A{\mathbb{A}}
\def\Af{\A^{\infty}}
\def\rhobarproj{\widetilde{\rho}}
\def\m{\mathfrak{m}}
\def\H{{\mathbb H}^3}

\def\Ok{\mathcal O}
\def\Od{B}
\newcommand{\QQ}{Q}
\def\A{\mathbb A}
\def\C{\mathbb C}
\def\F{\mathbb F}

\def\Q{\mathbb Q}
\def\R{\mathbb R}
\def\Z{\mathbb Z}

\def\Kbar{\osl{K}}
\def\Qbar{\ol{\Q}}
\def\Lbar{\ol{L}}

\def\Spec{\mathrm{Spec}}

\def\Frob{\mathrm{Frob}}
\newcommand{\maps}{\colon\thinspace}
\DeclareMathOperator{\vol}{Vol}
\DeclareMathOperator{\Isom}{Isom}
\DeclareMathOperator{\injrad}{injrad}

\newcommand{\setdef}[2]{{  \left\{  {#1}  \ \left| \   {#2} \right. \right\} }}
\newcommand{\abs}[1]{{\left| #1 \right|}}
\DeclareMathOperator{\tr}{tr}
\newcommand{\cusp}{\mathit{cusp}}

\newcommand{\RP}{{\R}{\mathrm{P}}}

\makeatletter
\def\cnewtheorem#1[#2]#3{\newtheorem{#1}{#3}[section]%
\expandafter\let\csname c@#1\endcsname\c@subsection}
\cnewtheorem{theorem}[subsection]{Theorem}
\cnewtheorem{conjecture}[subsection]{Conjecture}
\cnewtheorem{conj}[subsection]{Conjecture}
\cnewtheorem{lemma}[subsection]{Lemma}
\cnewtheorem{corollary}[subsection]{Corollary}
\cnewtheorem{cor}[subsection]{Corollary}
\cnewtheorem{proposition}[subsection]{Proposition}
\cnewtheorem{prop}[subsection]{Proposition}
\cnewtheorem{claim}[subsection]{Claim}
\cnewtheorem{question}[subsection]{Question}
\cnewtheorem{criterion}[subsection]{Criterion}
\cnewtheorem{goal}[subsection]{Goal}

\cnewtheorem{virtual_haken}[subsection]{Virtual Haken Conjecture}
\cnewtheorem{vpbn}[subsection]{Virtual Positive Betti Number Conjecture}
\makeautorefname{virtual_haken}{Conjecture}
\makeautorefname{vpbn}{Conjecture}

\theoremstyle{remark}
\cnewtheorem{remark}[subsection]{Remark}

\theoremstyle{definition}
\cnewtheorem{definition}[subsection]{Definition}

  \let\c@equation\c@subsection
  \let\p@equation\p@subsection
  \def\theequation{\bf \thesubsection}
  \let\cl@equation\cl@subsection
\makeatother

\begin{document}

\begin{asciiabstract}
We investigate a question of Cooper adjacent to the Virtual Haken
Conjecture.  Assuming certain conjectures in number theory, we show that
there exist hyperbolic rational homology 3-spheres with arbitrarily
large injectivity radius. These examples come from a tower of abelian
covers of an explicit arithmetic 3-manifold. The conjectures we must
assume are the Generalized Riemann Hypothesis and a mild strengthening
of results of Taylor et al  on part of the Langlands Program for GL_2
of an imaginary quadratic field.

The proof of this theorem involves ruling out the existence of
an irreducible two dimensional Galois representation (rho) of
Gal(Qbar/Q(sqrt(-2))) satisfying certain prescribed conditions. In
contrast to similar questions of this form, (rho) is allowed to have
arbitrary ramification at some prime of Z[sqrt(-2)].

Finally, we investigate the congruence covers of twist-knot orbifolds. Our
experimental evidence suggests that these topologically similar orbifolds
have rather different behavior depending on whether or not they are
arithmetic. In particular, the congruence covers of the nonarithmetic
orbifolds have a paucity of homology.
\end{asciiabstract}

\begin{webabstract} 
We investigate a question of Cooper adjacent to the Virtual Haken
Conjecture.  Assuming certain conjectures in number theory, we show
that there exist hyperbolic rational homology 3--spheres with
arbitrarily large injectivity radius.  These examples come from a
tower of abelian covers of an explicit arithmetic 3--manifold.  The
conjectures we must assume are the Generalized Riemann Hypothesis
and a mild strengthening of results of Taylor et al on part of the
Langlands Program for $\mathrm{GL}_2$ of an imaginary quadratic field.

The proof of this theorem involves ruling out the existence of an
irreducible two dimensional Galois representation $\rho$ of
$\mathrm{Gal}(\overline{\mathbb{Q}}/\mathbb{Q}(\sqrt{-2}))$ satisfying
certain prescribed ramification conditions.  In contrast to similar
questions of this form, $\rho$ is allowed to have arbitrary ramification
at some prime $\pi$ of $\mathbb{Z}[\sqrt{-2}]$.

In the next paper in this volume, Boston and Ellenberg apply pro--$\!p$
techniques to our examples and show that our result is true
unconditionally.  Here, we give additional examples where their
techniques apply, including some non-arithmetic examples.

Finally, we investigate the congruence covers of twist-knot
orbifolds.  Our experimental evidence suggests that these
topologically similar orbifolds have rather different behavior
depending on whether or not they are arithmetic.  In particular, the
congruence covers of the non-arithmetic orbifolds have a paucity of
homology.
\end{webabstract}

\begin{htmlabstract} 
<p class="noindent">
We investigate a question of Cooper adjacent to the Virtual Haken
Conjecture.  Assuming certain conjectures in number theory, we show
that there exist hyperbolic rational homology 3&ndash;spheres with
arbitrarily large injectivity radius.  These examples come from a
tower of abelian covers of an explicit arithmetic 3&ndash;manifold.  The
conjectures we must assume are the Generalized Riemann Hypothesis
and a mild strengthening of results of Taylor et al on part of the
Langlands Program for GL<sub>2</sub> of an imaginary quadratic field.</p>

<p class="noindent">
The proof of this theorem involves ruling out the existence of an
irreducible two dimensional Galois representation &rho; of
Gal(<span
style="text-decoration:overline"><b>Q</b></span>/<b>Q</b>(&radic;(-2)))
satisfying certain prescribed
ramification conditions.  In contrast to similar questions of this
form, &rho; is allowed to have arbitrary ramification at some prime
&pi; of <b>Z</b>[&radic;(-2)].</p>

<p class="noindent">
In the next paper in this volume, Boston and Ellenberg apply pro&ndash;p
techniques to our examples and show that our result is true
unconditionally.  Here, we give additional examples where their
techniques apply, including some non-arithmetic examples.</p>

<p class="noindent">
Finally, we investigate the congruence covers of twist-knot
orbifolds.  Our experimental evidence suggests that these
topologically similar orbifolds have rather different behavior
depending on whether or not they are arithmetic.  In particular, the
congruence covers of the non-arithmetic orbifolds have a paucity of
homology.</p>
\end{htmlabstract}

\begin{abstract} 
We investigate a question of Cooper adjacent to the Virtual Haken
Conjecture.  Assuming certain conjectures in number theory, we show
that there exist hyperbolic rational homology 3--spheres with
arbitrarily large injectivity radius.  These examples come from a
tower of abelian covers of an explicit arithmetic 3--manifold.  The
conjectures we must assume are the Generalized Riemann Hypothesis
and a mild strengthening of results of Taylor et al on part of the
Langlands Program for $\GL_2$ of an imaginary quadratic field.

The proof of this theorem involves ruling out the existence of an
irreducible two dimensional Galois representation $\rho$ of
$\Gal\bigl(\Qbar/\Q\bigl(\sqrt{-2}\bigr)\bigr)$ satisfying certain prescribed
ramification conditions.  In contrast to similar questions of this
form, $\rho$ is allowed to have arbitrary ramification at some prime
$\pi$ of $\Z[\sqrt{-2}]$.

In the next paper in this volume, Boston and Ellenberg apply pro--$\!p$
techniques to our examples and show that our result is true
unconditionally.  Here, we give additional examples where their
techniques apply, including some non-arithmetic examples.

Finally, we investigate the congruence covers of twist-knot
orbifolds.  Our experimental evidence suggests that these
topologically similar orbifolds have rather different behavior
depending on whether or not they are arithmetic.  In particular, the
congruence covers of the non-arithmetic orbifolds have a paucity of
homology.
\end{abstract}

\maketitle


\section{Introduction}

Here we investigate questions adjacent to the Virtual Haken
Conjecture, which concerns  (immersed) surfaces in 3--manifolds.  Let
$M$ be a closed 3--manifold, that is, one that is compact and has no
boundary.  An embedded orientable surface $S \neq S^2$ in $M$ is
\emph{incompressible} if $\pi_1(S) \to \pi_1(M)$ is injective.  The
manifold $M$ is called \emph{Haken} if it is irreducible and contains
an incompressible surface.  A Haken 3--manifold necessarily has
infinite fundamental group, but there are many such manifolds which
are not Haken.  One of the most interesting conjectures about
3--manifolds is Waldhausen's Conjecture \cite{Waldhausen68} which
posits the following:
\begin{virtual_haken}\setobjecttype{Conj}\label{conj:haken}
  Suppose $M$ is an irreducible 3--manifold with infinite fundamental
  group.  Then $M$ has a finite cover which is Haken.  
\end{virtual_haken}
Here, the term ``virtual'' refers to being allowed to pass to finite
covers.  For such questions, we can always take $M$ to be orientable,
and will do so from now on.  Assuming the Geometrization Conjecture, a
proof of which has been announced by Perelman \cite{Perelman1,Perelman2}, the unknown (and generic!) case of
\fullref{conj:haken} is when $M$ is hyperbolic, that is, has a
Riemannian metric of constant curvature $-1$.  Equivalently, $M = \H
/ \Gamma$ where $\H$ is hyperbolic 3--space, and $\Gamma$ is a
torsion-free uniform lattice in $\Isom^+(\H) \cong \PSL_2(\C) \cong \PGL_2(\C)$.

We focus on a stronger form of \fullref{conj:haken}, which
posits the existence of interesting homology in a finite cover:
\begin{vpbn}\setobjecttype{Conj}\label{vpbn}
  Let $M$ be a closed hyperbolic 3--manifold.  Then $M$ has a finite
  cover $N$ where the Betti number
  $$\beta_1(N) = \dim H^1(N ; \Q)> 0.$$
\end{vpbn}
The connection to the original conjecture is that $\beta_1(N) > 0$ implies
by Poincar\'e duality that $H_2(N; \Z) \neq 0$, and any non-trivial
element of the latter group can be represented by an incompressible
surface.  There are now many classes of examples for which
\fullref{vpbn} is known to hold, but there seems to be no
general approach.  In the case that $\Gamma = \pi_1(M)$ is an
\emph{arithmetic} lattice, \fullref{vpbn} can be naturally related
to automorphic forms.
   However, even in the arithmetic setting,
\fullref{vpbn} is known only in special cases, eg when the
field of definition has a subfield of index 2
(see Millson \cite{Millson}, Labesse--Schwermer \cite{LabesseSchwermer},
Clozel \cite{Clozel}, Lubotzky \cite{LubotzkyArith}, and Rajan
\cite{Rajan}).

In Kirby's problem list \cite[Problem~3.58]{Kirby}, Cooper formulated
the following question, which we devote this paper to studying.  The
\emph{injectivity radius} of $M$, denoted $\injrad(M)$, is the radius of
the largest ball that can be embedded around \emph{every} point in $M$;
equivalently, it is half the length of the shortest closed geodesic in
$M$.
\begin{question}[Cooper]\label{question:cooper}
  Does there exist a constant $K$ such that if $M$ is a closed hyperbolic
  3--manifold with $\injrad(M) > K$, then $\beta_1(M) > 0$?
\end{question}
A yes answer to this question would imply \fullref{vpbn},
because any hyperbolic 3--manifold has a congruence cover with
injectivity radius bigger than a fixed $K$.  However, the general
expectation was that the answer to this question is no; Cooper's
motivation in formulating it was to illustrate the depth of our
ignorance about \fullref{vpbn}.  The point of this paper is to
show that the answer to \fullref{question:cooper} is indeed no,
assuming certain conjectures in number theory.  In particular, we need
to assume the Generalized Riemann Hypothesis (GRH) and a mild
strengthening of results of Taylor et al on part of the Langlands
Program (\fullref{conj:global}).  

A \emph{rational homology
  sphere} is a closed orientable 3--manifold with $\beta_1 = 0$.  The
topological side of our main result is:
\begin{theorem}\label{thm:main-top}
  Assume the GRH and \fullref{conj:global}.  Then there exists
  an explicit tower of covers
  \[
  N_0 \gets N_1 \gets N_2 \gets N_3 \gets \cdots 
  \]
  so that each $N_n$ is a hyperbolic rational homology sphere, and
  $\injrad(N_n) \to \infty$ as $n \to \infty$.  Each cover $N_n \gets N_{n+1}$ is a
  regular cover with covering group either $\Z/3\Z$ or $\Z/3\Z \oplus
  \Z/3\Z$.  Moreover, the composite cover $N_0 \gets N_n$ is regular.
\end{theorem}
In the language of lattices, we have a nested sequence of lattices
$\Gamma_n = \pi_1(N_n)$ in $\PSL_2(\C)$ such that
$\bigcap_{n=1}^\infty \Gamma_n = 1$ and $H^1(\Gamma_n; \R) = 0$ for
each $n$ (equivalently, the abelianization of each $\Gamma_n$ is
finite).  Examples of infinite towers of covers of hyperbolic rational
homology spheres have been previously constructed by Baker, Boileau,
and Wang \cite{BakerEtAl}, but these lack the crucial requirement on
the injectivity radius.

Despite \fullref{thm:main-top}, the
$N_n$ are easily seen to satisfy \fullref{vpbn}.  Amusingly, it
turns out that all of the $N_n$ are in fact Haken (see
\fullref{subsec:virt-prop}).  Thus they do not give a no answer to
\begin{question}
  Does there exist a constant $K$ such that if $M$ is a closed hyperbolic
  3--manifold with $\injrad(M) > K$, then $M$ is Haken?
\end{question}
As with Cooper's original question, we suspect the answer must be  no, but
see no way of showing this. 

The construction of the examples in \fullref{thm:main-top} is 
 arithmetic in nature.   A precise statement is 
\begin{theorem}\label{thm:main-arith}
  Let $D$ be the \emph{(}unique\emph{)}
  quaternion algebra over $K = \Q(\sqrt{-2})$ ramified at $\pi$ and
  $\pibar$, where $3 = \pi \pibar$. Let $\Od$ be a maximal order of $D$.
  Let $\m$ be a maximal bi-ideal of $\Od$ trivial away from $\pi$.
  Finally, let $B_n$ be the complex embedding of $\Od^{\times} \cap (1 + \m^n)$ into
  $\PGL_2(\C)$, and let $M_n = \H/B_n$. Then assuming the Langlands
  conjecture for $\GL_2(\A_K)$ (\fullref{conj:global}) and the
  GRH, then $H^1(M_n ; \Q) = 0$ for all $n$.  
\end{theorem}
The examples $N_n$ of \fullref{thm:main-top} are the above $M_n$
with the first few dropped for technical reasons; see
\fullref{subsec:M_n} for details.  In the next paper in this
volume, Boston and Ellenberg apply pro--$\!p$ techniques to the above $M_n$
and show that $H^1(M_n ; \Q) = 0$ unconditionally
\cite{BostonEllenberg}.  We give additional examples where their
techniques apply, including some non-arithmetic examples, in
\fullref{sec-pro-p}.

It has long been known that the cohomology groups of arithmetic
3--manifolds are related to spaces of automorphic forms. This relationship
has been previously exploited: Clozel~\cite{Clozel} was able to prove
\emph{non-vanishing} results for certain cohomology groups by using
automorphic forms associated to  Gr\"{o}ssencharakters.  In contrast,
computations made by
 Grunewald, Helling, and Mennicke~\cite{Germans},
and Cremona~\cite{Cremona}, working with congruence covers of Bianchi
manifolds, suggest the paucity of automorphic forms, and lead to the
suspicion that for certain groups,  infinitely many congruence covers
contain no interesting (non-cuspidal) cohomology.  To this point such
problems have been unapproachable.
 The trace formula, so useful in many other situations,
here only contributes the equality $0 = 0$ (Poincare duality identifies
the two ``interesting'' cohomology groups, $H^1$ and $H^2$). In this
paper, we present an approach to these questions that succeeds (assuming
standard conjectures) in ruling out the existence of  automorphic forms
associated to particular congruence covers of arbitrarily large degree.

The starting point is the theorem of Taylor
et~al~\cite{taylorothers,taylorimag} which establishes
for the automorphic forms of interest a family of compatible
$\Gal(\Kbar/K)$--representations for some imaginary quadratic
field $K$.  The nature of these representations is not yet completely
understood, and in particular we must assume a slight strengthening of
their results. The main idea is that one can rule out the existence
of certain automorphic forms by ruling out the existence of the
corresponding Galois representations, even allowing the conductor of
these representations to become arbitrarily large.  This is in contrast
to the situation over $\Q$, where applications tend to work in reverse:
the proof of Fermat's last theorem \cite{TW,Wiles} uses the
 non-existence
of  modular forms of weight two and level $\Gamma_0(2)$
 to rule out the existence of
certain Galois representations! We restrict our representations $\rho$
in several steps. First, choose a prime $p = \pi \pibar$ that splits in
the ring of integers of $K$
 and consider automorphic forms  whose conductor is highly divisible
by $\pi$ but strictly controlled at $\pibar$;
 let $\rho$ be the associated   $\lambda$--adic
representation for some $\lambda | p$.  Second, using an idea of
Tate~\cite{Tate} we prove in a specific case that $\rhobar$ must be
reducible.  For this step we require the Generalized Riemann Hypothesis.
The key is now to use the local behavior of $\rho$ at $\pibar$ as a
fulcrum to prove that $\rho$ is reducible, which is enough to force a
contradiction. This last idea in a different guise can be seen in the
work of Fontaine~\cite{fontaine} (further developed by Schoof and others,
see for example Brumer--Kramer \cite{Brumer}, Calegari \cite{Semistable}
and Schoof \cite{schoof}).
%
%
%
%
Finally, our application to rational homology spheres comes
from switching between non-compact covers of Bianchi manifolds
and compact arithmetic quotients of $\H$ by using the
Jacquet--Langlands correspondence. We learnt of this idea
from Alan Reid.

Finally, we investigate the congruence covers of twist-knot orbifolds,
with the goal being to gather experimental evidence about how common
it is for congruence covers to have $\beta_1 > 0$.  For instance, one
would like to know if it is plausible to attack \fullref{vpbn}
solely by examining covers of this kind.  (For comparison, complex
hyperbolic manifolds are also expected to satisfy
\fullref{vpbn}.  However, there are infinitely many
incommensurable \emph{arithmetic} complex hyperbolic
manifolds all of whose congruence covers have $\beta_1 = 0$; see
Rapoport--Zink \cite{RapoportZink}, Rogawski \cite{Rogawski} and Clozel
\cite{Clozel2}.)  Our examination of a family of
topologically similar orbifolds found strikingly different behavior
depending on whether or not the orbifold was arithmetic.  In
particular, the congruence covers of the non-arithmetic orbifolds have
a paucity of homology.  Indeed, it seems plausible that our sample
includes a non-arithmetic hyperbolic orbifold where only finitely many
congruence covers of the form $\Gamma_0(\mathfrak{p})$ have $\beta_1 > 0$.

\subsection{Outline of contents}

\fullref{sec:examples} gives a detailed description of the
orbifolds $M_n$ and $N_n$ of the main theorems from both topological
and arithmetic points of view.  \fullref{sec:examples} also
derives \fullref{thm:main-top} from \fullref{thm:main-arith}.
\fullref{sec:modular-forms} discusses the connection to
automorphic forms, states the needed form of the Langlands conjecture,
and gives the reduction of \fullref{thm:main-arith} to showing the
non-existence of certain Galois representations.  Such Galois
representations are ruled out in
Sections~\ref{sec:Serre-level}--\ref{sec:global}.  An elaboration of
the Boston--Ellenberg pro--$\!p$ approach is given in
\fullref{sec-pro-p}, together with a list of examples to which it
applies.  Finally, \fullref{sec:twist} contains the experimental
results on the twist-knot orbifolds.

\subsection{Acknowledgments}  

Calegari was partially supported as a 5--year fellow of the American
Institute of Mathematics.  Dunfield was partially supported by U.S.
National Science Foundation grant \#DMS-0405491, as well as a Sloan
Fellowship.  The authors thank Nigel Boston, Kevin Buzzard, Jordan
Ellenberg, Oliver Goodman, Damian Heard, Craig Hodgson, Alex Lubotzky,
Hee Oh, Dinakar Ramakrishnan, Alan Reid, and Richard Taylor for
helpful conversations and correspondence.

\section{The examples}\label{sec:examples}

In this section, we give a detailed description of the orbifolds $M_n$
of \fullref{thm:main-arith}, and also derive
\fullref{thm:main-top} from \fullref{thm:main-arith}.  

\subsection{The orbifold $M_0$}\label{subsec:M0}

We start by looking at $M_0$, which is constructed arithmetically as
follows (for details see Machlachlan--Reid \cite{MacReid}, and Vigneras
\cite{Vigneras}).  We start
with the field $K = \Q(\sqrt{-2})$.  The prime $3$ splits over $K$
into two primes $\pi = 1 - \sqrt{-2}$ and $\pibar = 1 + \sqrt{-2}$.
Let $D$ be the unique quaternion algebra over $K$ which is ramified at
exactly the places $\pi$ and $\pibar$.  Explicitly, $D$ has a standard
basis $\{ 1, i, j, \iijj \}$ where $i^2 = -1, j^2 = -3$, and $\iijj = -ji$;
equivalently, $D$ is said to be given by the Hilbert symbol $\left(
  \frac{-1, -3}{K} \right)$.  Let $\Od$ be a maximal order in $D$.
The order $\Od$ is unique up to conjugacy, as the number of conjugacy
classes of maximal orders divides the restricted class number $h_\infty$
of $K$, which here is just the class number $h = 1$ (see
eg Machlachlan--Reid \cite[Section~6.7]{MacReid}).  Now let $\Od^\times$ be the group of units of
$\Od$.  Taking the complex place of $K$ gives us an embedding $D \to D
\otimes_K \C \cong M_2(\C)$ into the algebra of $2 \times 2$ matrices over $\C$.
The image of $\Od^\times$ lands in $\GL_2(\C)$, and composing this with the
projectivization $\GL_2(\C) \to \PGL_2\C$ gives a homomorphism $\rho
\maps \Od^\times \to \PGL(2,\C)$.  The target group $\PGL_2\C$ can naturally
be identified with the group of orientation preserving isometries of
hyperbolic 3--space $\H$.  The group $\rho(\Od^\times)$ is a lattice, and the
base orbifold for our examples is
\[
M_0 = \H / \rho(\Od^\times).
\]
Since, by construction, $D$ ramifies at some finite place, the orbifold $M_0$
is compact.  

\begin{remark} 
  In the topology literature, it is more typical to consider $\Od^1$,
  the elements of $\Od$ of reduced norm $1$, rather than $\Od^\times$. At the
  complex place, such elements lie in $\SL_2(\C) \leq \GL_2(\C)$, and 
  also give a lattice $\rho(\Od^1) \leq \rho(\Od^\times)$ in $\PGL_2(\C)$.  The
  reduced norm gives a homomorphism from $\Od^\times$ to the group of units
  of $K$, which is just $\pm1$.  The group $\Od^1$ is the kernel, and
  hence $\left[ \rho(\Od^\times) : \rho(\Od^1)\right] = 2$.  Thus, if we set
  $M_0' = \H / \rho(\Od^1)$, we that have $M_0'$ is a 2--fold cover
  of $M_0$.  From the point of view of automorphic forms, it is
  more natural to work with $M_0$ rather than $M_0'$.
\end{remark}

\begin{figure}[ht!]
\labellist\small
\pinlabel {$4$} [r] at 0 107
\pinlabel {$4$} [rb] at 64 178
\pinlabel {$2$} [l] at 283 107
\pinlabel {$2$} [lb] at 133 140
\pinlabel {$3$} [t] at 172 102
\endlabellist
\centerline{\includegraphics[scale=0.63]{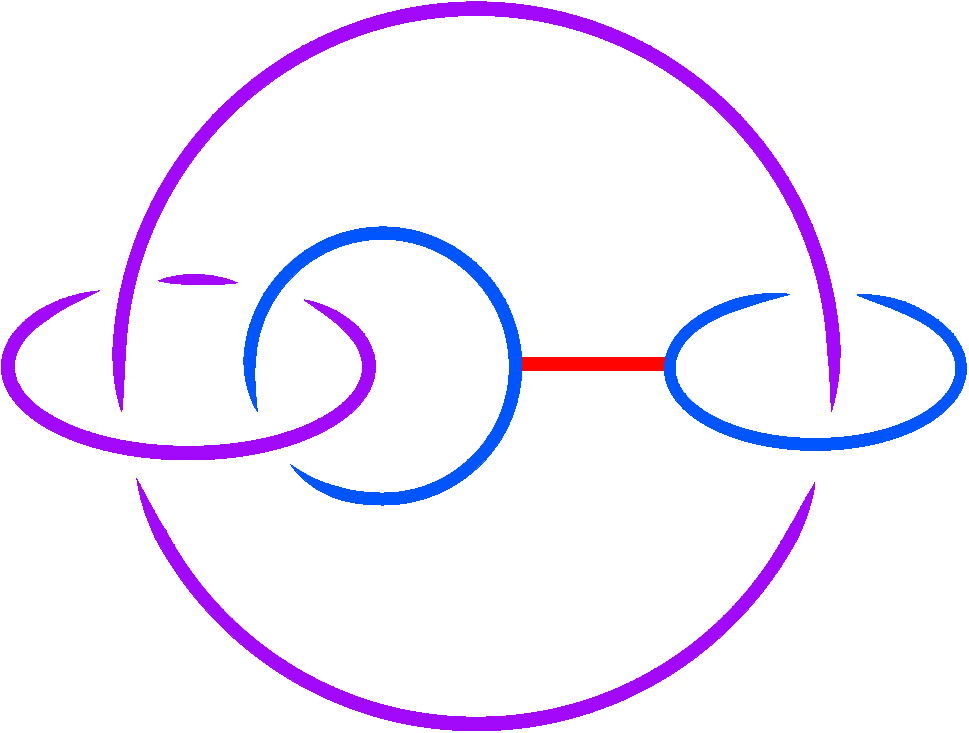}}
\caption{The orbifold $M_0 = \H /  \rho(\Od^\times)$ has underlying
space $S^3$ with orbifold locus the above graph.  The labels on the
edges of the graph correspond to the order of the cyclic group action
associated to that part of the singularity.}
\label{fig:pgl-orbifold}
\end{figure}
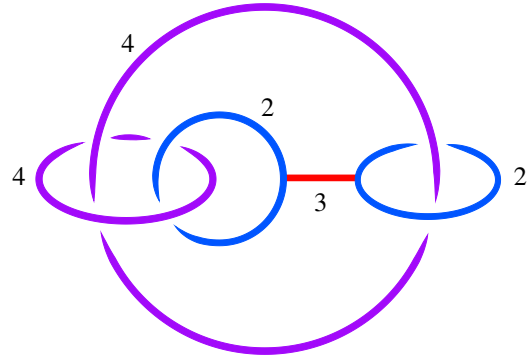

The basic topology of $M_0$ is also easy to describe.  The underlying
topological space for $M_0$ is just the 3--sphere $S^3$, and the
orbifold locus is shown in \fullref{fig:pgl-orbifold} (for the
derivation of this from the arithmetic, see
\fullref{subsec:top-der}).  The hyperbolic volume of $M_0$ can
be directly calculated from data about the quaternion algebra, see
eg Machlachlan--Reid \cite[Theorem~11.1.3]{MacReid}:
\begin{equation}\label{eqn:volume}
\vol(M_0) = \frac{1}{2} \vol(M_0') = \frac{8 \sqrt{2}}{\pi^2} \zeta_K(2) = 2.0076820066823962745447297\ldots
\end{equation}
Here is a presentation for $\pi_1(M_0) = \rho(\Od^\times)$:
\begin{equation}\label{eq:presentation}
\bigl\langle u,v,x,y ~\big|~ u^2,v^2, x^4, y^4, y x y^{-1}  v x^{-1} v,
x^{-1}  v  x  v  u  y^{-1}  u  y, (u y^{-1} u y)^3 \bigr\rangle
\end{equation}
which is obtained by eliminating excess generators from a Wirtinger
presentation derived from \fullref{fig:pgl-orbifold}. Let's now
give the canonical representation $\pi_1(M_0) \to \PSL_2(\C)$ in terms
of the quaternion algebra.  The maximal order $\Od$ can be taken to be
the $\Ok_K$--span of 
\[
\bigl\langle 1, i, s=(1/2)(i + j), t = (1/2)(1 +  \iijj)\bigr\rangle.
\]  
As discussed in \fullref{subsec:top-der}, in terms of our
presentation for $\pi_1(M_0)$ the unit group  $\Od^\times$ is generated by
\begin{align*}
u &\mapsto i, & v &\mapsto -2 i + \pi  s, \\
x &\mapsto \sqrt{-2} \cdot  1  + \pi i  + \pibar s  - \sqrt{-2} t, &  y &\mapsto \sqrt{-2} \cdot 1 + s - \sqrt{-2} t.
\end{align*}
The representation $\pi_1(M_0) \to \PSL_2(\C)$ is then obtained by taking
an explicit embedding of $D$ into $M_2(\C)$, eg
\[
i \mapsto \begin{pmatrix} i & 0 \\ 0 & -i  \end{pmatrix}, \quad %
j \mapsto \begin{pmatrix}  0 & 1 \\ -3 & 0 \end{pmatrix}.
\]

\subsection{The orbifolds $M_n$}\label{subsec:M_n}

We turn now to the orbifolds $M_n$.  By definition $M_n$ is the
congruence cover of $M_0$ of level $\pi^n$.  Very succinctly, this
means that if $\m$ is the maximal bi-ideal of $\Od$ trivial away from
$\pi$, then $M_n = \H / \rho\big( \Od^\times \cap (1 + \m^n) \big)$.  We will
now describe these congruence covers in more detail, so that we may
derive \fullref{thm:main-top} from \fullref{thm:main-arith}.
It is worth noting that since these are congruence covers for a prime
that ramifies in $D$, their structure is rather different than those for
a generic prime.

As usual, $K_\pi$ will denote the completion of $K$ with respect to the
$\pi$--adic norm.  The integers in $K$ will be denoted $\Ok =
\Z[\sqrt{-2}]$, the $\pi$--adic completion of $\Ok$
is  the valuation ring $\Ok_\pi$.
As $\Ok$ has unique factorization, we conflate primes ideals and prime
elements, so that $\pi$ is also a uniformizing element in $\Ok_\pi$, with $\pi
\Ok_\pi$ the unique maximal ideal.  Now consider $D$ at this finite place,
which we denote $D_\pi = K_\pi \otimes_K D$.  That $D$ ramifies at $\pi$ means
precisely that $D_\pi$ is the unique division algebra over $K_\pi \cong \Q_3$.
Explicitly, $D_\pi$ is the following (for details, see Maclachlan--Reid
\cite[Section~2.6, Section~6.4]{MacReid} which we follow closely).  Let $L$ be the unique
unramified quadratic extension of $K_\pi$; then $L = K_\pi(\sqrt u)$ for
some unit $u \in \Ok_\pi^\times$.  Explicitly, $D_\pi$ is specified by the
Hilbert symbol $\bigl( \frac{u,\pi}{K_\pi} \bigr)$.  Even more
concretely, one can take
\begin{equation}\label{eq:Dpi-concrete}
D_\pi =  \biggl\{ \,  \begin{pmatrix} a & b \\ \pi b' & a' \end{pmatrix}
\, \biggr\} \quad \mbox{where $ a, b \in L$}
\end{equation}
  and $'$ denotes the Galois involution of  $L / K$.
In this model, $\{1, i, j \}$ correspond to $(a,b) = (1, 0),
(\sqrt{u}, 0)$, and $(0, 1)$ respectively.  A concise way of
writing an element of $D_\pi$ is thus $a + b j$ where $a,b \in L$.
The reduced norm $n$ is then $aa' - \pi bb'$, and the trace $a +
a'$.

The algebra $D_\pi$ has a natural norm $\abs{d} = \abs{n(d)}_\pi$ where
$\abs{\, \cdot \,}_\pi$ is the norm on $K_\pi$.  The valuation ring $\Od_\pi
= \big\{d \in D_\pi \: \big\vert \; \abs{d} \leq 1\big\}$ is the unique
maximal order of $D_\pi$.  In \eqref{eq:Dpi-concrete}, the elements of
$\Od_\pi$ are simply those that have $a,b \in \Ok_L$.  We have a natural
embedding $\Od{} \hookrightarrow \Od_\pi$; as $\Od$ is maximal in $D$, the valuation
ring $\Od_\pi$ is equal to $\Od{} \otimes_\Ok \Ok_\pi$.  The unique maximal
bi-ideal of $\Od_\pi$ is just $\QQ = \Od_\pi j$.  The units $\Od^\times_\pi$
have a natural filtration
\[
\Od^\times_\pi \supset 1 + \QQ \supset 1 + \QQ^2 \supset 1 + \QQ^3 \supset  \cdots 
\]
Let $\Gamma_n$ be the preimage of $1 + \QQ^n$ under $\Od^\times \hookrightarrow \Od_\pi$, and set
$M_n = \H / \rho(\Gamma_n)$.  By definition, $M_n$ is the congruence
cover of $M_0$ of level $\pi^n$.  To relate this back to language at
the beginning of this subsection, the bi-ideal $\m$ of $\Od$ is
precisely the preimage of $\QQ$ under $\Od{} \hookrightarrow \Od_\pi$, and so $\Gamma_n = \Od^\times
\cap (1 + \m^n)$.

We now examine the orbifolds $M_n$ more closely, and so derive
\fullref{thm:main-top} from \fullref{thm:main-arith}.
\begin{proof}[Proof of \fullref{thm:main-top}]
  The manifold $N_n$ in the statement of the theorem is simply
  $M_{n+d}$ for a fixed positive integer $d$.  (We need to drop the
   first few $M_n$ as they are genuine orbifolds, not manifolds.)  Thus
  modulo \fullref{thm:main-arith}, we need to check three things:
  \begin{enumerate}
  \item The quotient $M_n = \H / \rho(\Gamma_n)$ is a manifold for large
    $n$; equivalently $\rho(\Gamma_n)$ is torsion free.

  \item The injectivity radius of $M_n$ goes to $\infty$ as $n \to \infty$.
    
  \item The covering group of $M_{n} \gets M_{n + 1}$ is $\Z/3\Z$ for
    even $n > 0$ and $(\Z/3\Z)^2$ for odd $n > 0$.  Moreover, $M_{0} \gets
    M_{n}$ is a regular cover for all $n$.

  \end{enumerate}
  First, let us examine $M_0 \gets M_1$.  Consider the residue field of
  $L$, denoted $\Lbar = \Ok_L/\pi \Ok_L \cong \F_9$.   Then 
  \[
  \Od^\times_\pi / (1 + \QQ) \cong \Lbar^\times \cong \Z/8\Z %
  \quad \mbox{via $a + b j \mapsto a + \pi \Ok_L$}.  
  \]
  To determine $\left[ \Od^\times : \Gamma_1\right]$, we need to understand the
  image of $\Od^\times$ in $\Od^\times_\pi$.  This image is not dense as all
  elements of $\Od^\times$ have norm in $\Ok^\times = \{\pm1\}$, whereas the set of norms
  elements is $\Od^\times_\pi$ is infinite.  However, if we let $\Od^{\pm1}_\pi$ consist of
  those elements of norm $\pm1$, then strong approximation (see
  Machlachlan--Reid \cite[Theorem~7.7.5]{MacReid}) shows that $\Od^\times$
  is dense in $\Od^{\pm1}_\pi$.  It is
  not hard to see that $\Od^{\pm 1}/(1 + \QQ)$ is still all of
  $\Lbar^\times$;
  for this and subsequent calculations, the reader may find it
  convenient to note $D_\pi \cong \bigl( \frac{-1,3}{\Q_3} \bigr)$.
  Thus it follows that $\Od^\times \to \Z/8\Z$ is surjective.  Thus $\left[
    \Od^\times : \Gamma_1\right] = 8$; however $\left[\rho(\Od^\times) :
    \rho(\Gamma_1)\right] = 4$ as the kernel of $\rho \maps \Od^\times \to
  \PGL_2(\C)$ is $\{ \pm1 \}$, which maps non-trivially under $\Od^\times \to
  \Z/8\Z$.  Thus $M_0 \gets M_1$ is a $4$--fold cyclic cover, and $\rho
  \maps \Gamma_n \to \PGL_2{\C}$ is injective for all $n > 0$.  Notice also
  that the reduced norm map $n \maps \Od^\times \to \{ \pm 1 \}$ is the same
  as the composite $\Od^\times \to \Z/8\Z \to \Z/2\Z$; hence $\Gamma_n$ lies in
  $\Od^1$ for $n > 0$.

  Now look at $M_{n} \gets M_{n + 1}$.   In this case, we have
  \[
1 \to \Gamma_{n+1} \to   \Gamma_n \to \left(1 + \QQ^n\right) / \bigl(1 +
\QQ^{n + 1}\bigr) \cong \Lbar^+ \cong (\Z/3\Z)^2,
  \]
  where the interesting identification is given by $1 + (a + b j)j^n
  \mapsto a + \pi\Ok_L$.  The image of the rightmost map is $\Od_\pi^1/(1 +
  \QQ^{n+1})$, which turns out to be $\Z/3\Z$ when $n$ is even and
  $(\Z/3\Z)^2$ when $n$ is odd.  Thus $M_{n} \gets M_{n + 1}$ is has
  covering group either $\Z/3\Z$ or $(\Z/3\Z)^2$.  Also, since $\QQ^n$
  is a bi-ideal, it is immediate that $\Gamma_n$ is normal in $\Od^\times$.
  This establishes point (3) above.
  
  Turning now to (1) and (2), consider some $\Gamma_n$ for a fixed $n >
  0$.  As noted above, $\Gamma_n$ is in $\Od^1$.  Hence for $g \in \Gamma_n$, the
  trace of $g$ is the same the trace of $\rho(g)$ in $\PSL_2(\C)$ (the
  latter of which is only defined up to sign).  Suppose $g$ is
  non-trivial.  If $g$ has finite order, then $\tr (g)$ is a real
  number in $(-2, 2)$.  Otherwise $g$ corresponds to a closed geodesic in
  $M_n$.  If the length of this geodesic is $l$ and the twist parameter is $\theta \in [0, 2
  \pi)$, then
  \begin{equation}
  \cosh\left( \frac{l + i \theta}{2}\right) = \pm \frac{\tr(g)}{2}.
  \end{equation}
  Thus to prove (1) and (2) it suffices to show that 
  \begin{equation}\label{eq:trace-grows}
  \min\big\{\,  \abs{\tr(g)} \; \; \big| \; \; {g \in \Gamma_n \setminus \{1\}} \big\} \to \infty \quad \mbox{as $n \to \infty$,}
  \end{equation}
  where $\abs{\, \cdot \,}$ denotes the absolute value on $K$ at the complex place.
  Fix $g \in \Gamma_n \setminus \{1\}$.  Now, as $g = 1 + x j^n$, we have $\tr(g)
  = 2 + y \pi^m$ where $m = \left\lceil n/2\right\rceil$ and $y \in \Ok$.  Hence 
  \[
  \abs{\tr(g)} \geq \abs{y} \cdot \abs{\pi}^m - 2 \geq 1 \cdot 3^{m/2} - 2 \geq 3^{n/4} - 2, 
 \]
 and thus \eqref{eq:trace-grows} holds, completing the proof of the theorem.
 \end{proof}

\subsection{Finding a topological description of $M_0$}\label{subsec:top-der}

In this subsection, we outline the procedure used to find the
topological description of $M_0$ given in
\fullref{fig:pgl-orbifold}.  Along the way, we find topological
descriptions of $M_0' = \H / \rho(\Od^1)$ and $M_1$, the latter of
which has the interesting consequence that the examples of
\fullref{thm:main-top} are Haken (see
\fullref{subsec:virt-prop}). 

\begin{figure}[ht!]
\labellist\small
\pinlabel {$3$} [b] at 160 296
\pinlabel {$3$} [b] at 160 436
\pinlabel {$(4,1)$} [l] at 320 361
\pinlabel {$(2,1)$} [r] at 72 361
\pinlabel {$(2,1)$} [bl] at 280 420
\endlabellist
\centerline{\includegraphics[scale=0.90]{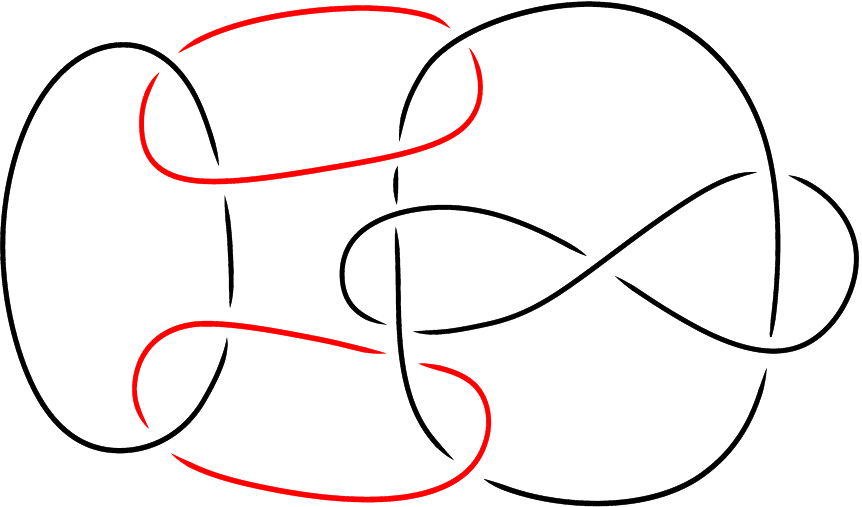}
}
\caption{Dehn surgery description of an orbifold $N$, which will be
shown to be $M_1$. Our framing conventions for Dehn surgery are given
    by: \protect\raisebox{-2pt}{\protect\reallyincludegraphics[scale=0.35,
      angle=90]{\figdir/orientation}}} 
  \label{fig:link1}
\end{figure}

When trying to determine the topology of $M_0$, it is important to
remember that the structure of a commensurability class of arithmetic
3--manifolds is quite complicated; in particular, it contains
infinitely many minimal orbifolds (throughout this subsection, see
Maclachlan--Reid \cite{MacReid} for details).  We started by simply finding some
arithmetic orbifold commensurable to $M_0$.  A natural place to look
is the Hodgson--Weeks census of small volume closed hyperbolic
3--manifolds \cite{SnapPea}; however, the smallest manifold
commensurable with $M_0$ has much too large a volume to appear there.
Instead, we started with cusped manifolds in the
Callahan--Hildebrand--Weeks census and did \emph{orbifold} Dehn filling,
looking for something with volume a rational multiple of the value given by
\eqref{eqn:volume}.  The commensurability class of an arithmetic
hyperbolic 3--orbifolds is completely determined by the invariant
trace field and quaternion algebra, which may be computed using
Goodman's program Snap \cite{Snap,SnapArticle}.  One thus finds that
the orbifold $s594(3,-3)$ is commensurable with $M_0$ and
$\vol(s594(3,-3)) = 2 \vol(M_0)$.  However, the (non-invariant) trace
field of $s594(3,-3)$ is bigger than $K$, so $s594(-3,3)$ is not
derived from a quaternion algebra; that is, it is not conjugate into
$\rho(\Od^1)$.  Passing to a 2--fold cover gives an orbifold $N$ which lies
in $\rho(\Od^1)$ with index 2.  Using Thistlethwaite's table of links
provided with \cite{Snap}, a brute force search finds the Dehn surgery
description of $N$ shown in \fullref{fig:link1}.

\begin{figure}[ht!]
\labellist\small
\pinlabel {$\pi$} [rb] at 93 361
\endlabellist
\centerline{\includegraphics[scale=0.72]{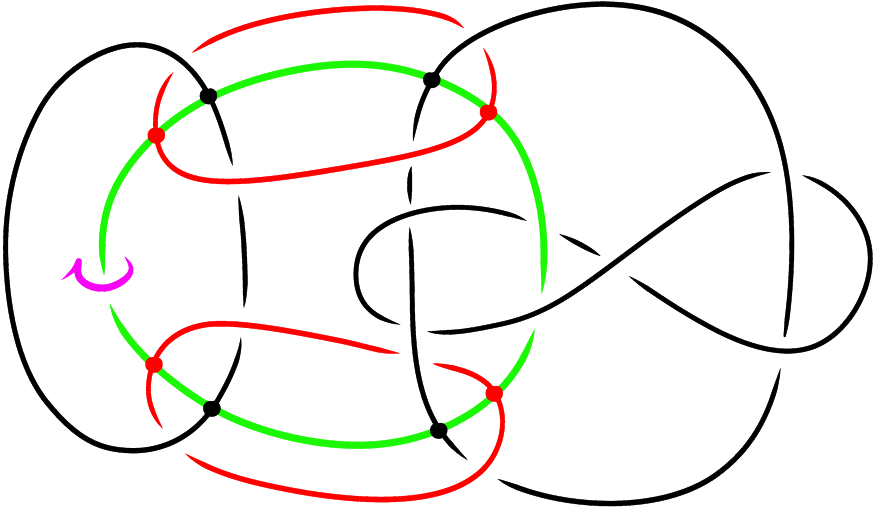}}
\caption{The involution $\tau$ which quotients $N$ down to $M_0'$.  The
fixed point set of $\tau$ intersects the link in $8$ points.}
\label{fig:link2}
\end{figure}

\begin{figure}[ht!]
\labellist\small
\pinlabel {$3$} [b] at 285 745
\pinlabel {$3$} [rt] at 285 570
\pinlabel {$(4,2)$} [lt] at 445 620
\endlabellist
\centerline{\includegraphics[scale=0.42]{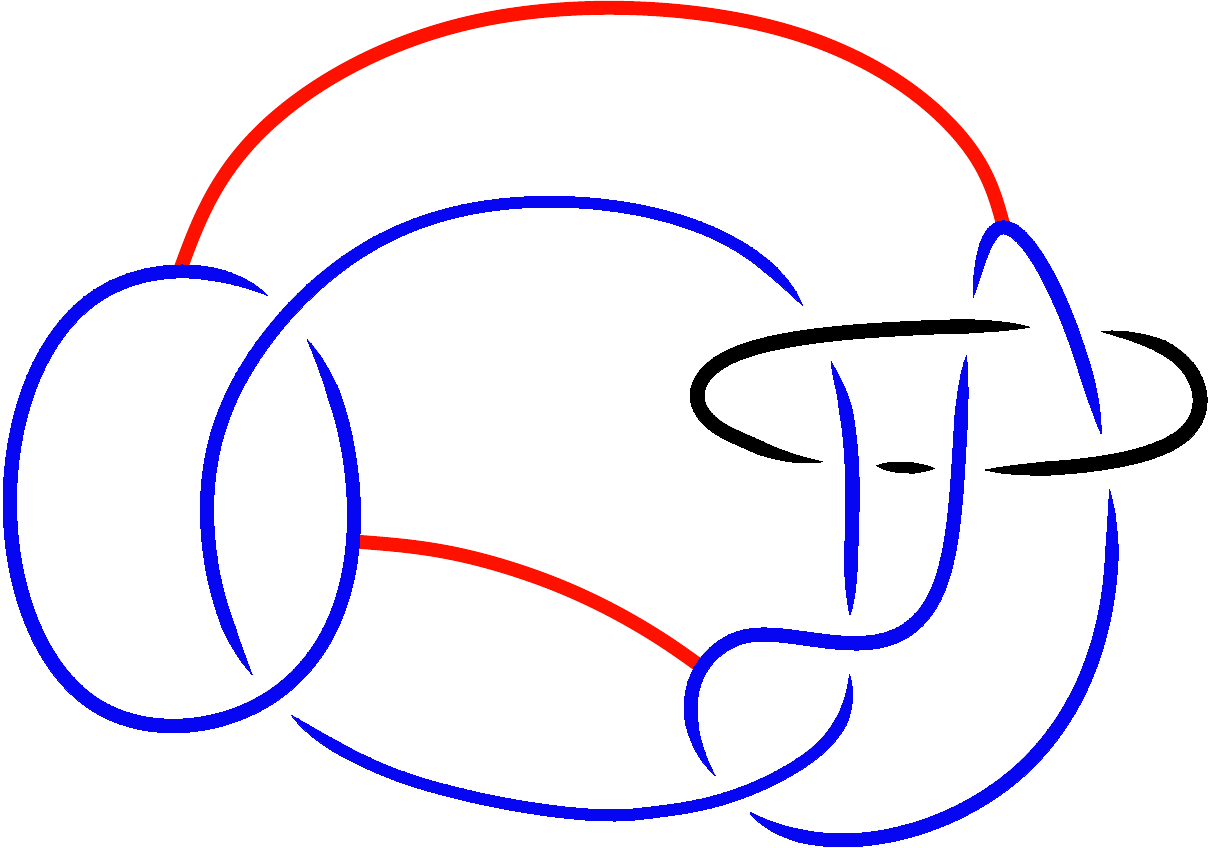}}
\caption{The orbifold $M_0' = \H / \rho(\Od^1)$.  The unlabeled part of
the singular graph should be labeled 2.}
\label{fig:psl-orbifold}
\end{figure}

Since $N$ is derived from a quaternion algebra, and $\vol(N) = 2
\vol(M_0')$, we know that there is some involution of $N$ which
quotients it down to $M_0'$.  Unfortunately, the symmetry group of $N$
is quiet large ($\Isom(N) = \Z/2\Z \times D_8$), and there are some $11$ distinct
involutions of $N$.  In order to find the correct one, we used Snap to
give matrices in $\PSL_2\C$ inducing each of these automorphisms of
$\pi_1(N) \leq \PSL_2\C$.  Adjoining these to $\pi_1(N)$ one at a time,
one finds there is a unique one, call it $\tau$, where the trace field
remains just $K$.  Looking at the action of $\tau$ on short geodesics of $N$
shows that $\tau$ in fact comes from a symmetry of the link in the
surgery diagram of \fullref{fig:link1}.  The symmetry group of the
link is much smaller, namely just $\Z/2\Z \oplus \Z/2\Z$.  One also finds that
$\tau$ fixes each component of the orbifold locus of $N$ at exactly two
points.  This now forces $\tau$ to be the involution of $N$ given in
\fullref{fig:link2}.  The quotient $N/\tau$ is then $M_0' = \H /
\rho(\Od^1)$, which is given in \fullref{fig:psl-orbifold}.

To find $M_0$ itself, one again searches for an additional symmetry
$\tau'$ to add to $N$ so that $\left\langle \pi_1(N), \tau, \tau' \right\rangle$ has
trace field $\Q(\sqrt{-2}, i)$.  One finds that there is essentially
only possible choice for $\tau'$, and this has order $4$ with $(\tau')^2 =
\tau$.  It is not hard to see a symmetry of \fullref{fig:link1},
which quotients down to \fullref{fig:pgl-orbifold}.  Since the
orbifold in \fullref{fig:pgl-orbifold} has 4--torsion, is must be
the quotient of $N$ under the action of a cyclic group of order $4$.
However, there are two distinct $\Z/4\Z$ subgroups of $\Isom(N)$ which
contain $\tau$, so one must do another check to see that
\fullref{fig:pgl-orbifold} is really $M_0'$.

We started with the presentation \eqref{eq:presentation} for the
fundamental group of the orbifold give in
\fullref{fig:pgl-orbifold}.  Eliminating variables, it is possible
to solve for a representation $\rho_0 \maps \pi_1 \to \PSL_2\C$ which
\emph{appears} to be the canonical representation.  From this, one
derives the quaternion algebra picture given at the end of
\fullref{subsec:M0}.  Thus a posteriori, $\rho_0(\pi_1)$ is a
discrete group.  Using SnapPea \cite{SnapPea} to compute a Dirichlet
domain shows that $\rho_0(\pi_1)$ is cocompact with the correct volume,
and hence the claimed generators for $\Od^\times$ really do generate it.
One then observes that gluing up this Dirichlet domain according to
the face pairings gives \fullref{fig:pgl-orbifold}; since
residually finite groups are Hopfian it follows that $\rho_0$ is
faithful.  More simply, one can use Heard's new program Orb
\cite{Orb} to see that the canonical representation is as claimed.  

To conclude this section, we show that $N$ is none other than $M_1$.
From the proof of \fullref{thm:main-top}, we know that $M_0 \gets
M_1$ is cyclic cover of degree 4.  Also, the traces of elements of
$\Gamma_1 \leq \Od^1$ are congruent to $2$ mod $\pi$.  It follows that
$\pi_1(M_1) = \rho(\Gamma_1)$ contains no $2$--torsion.  Now $N$ and
the congruence cover of $M_1'$ corresponding to $\pibar$ are also
$4$--fold cyclic covers of $M_0$ with no $2$--torsion.  There are only
two candidate homomorphisms $\pi_1(M_0) \to \Z/4\Z$ whose kernels
contain no $2$--torsion.  As $\Od$ is maximal, strong approximation
implies that $M_1$ and $M_1'$ are distinct covers of $M_0$, Thus if
$N$ is not $M_1'$, it must be equal to $M_1$.  A quick check with
SnapPea shows that $\pi(N)$ contains elements whose traces are
congruent to $0$ mod $\pibar$, and hence $N \not\cong M_1'$.  (Note:
SnapPea only gives elements in $\PSL_2\C$, where the trace is defined
only up to sign.  This does not matter since we are testing whether the
trace is $0$ in $\F_3$.)  Thus $N \cong M_1$ as desired.

\subsection{Virtual properties of  $M_0$}\label{subsec:virt-prop}

Despite \fullref{thm:main-top}, the orbifold $M_0$ does satisfy
\fullref{vpbn}.  Indeed, since $K = \Q(\sqrt{-2})$ has a
subfield of index two, there are in fact congruence covers with $\beta_1
> 0$ (see Labesse--Schwermer \cite{LabesseSchwermer} and Lubotzky
\cite{LubotzkyArith}).  One way to see this
directly is to start from the fact that $M_0$ contains an immersed
totally geodesic surface.  This follows because $D$ is ramified at
exactly the two primes sitting over $3$ (see Maclachlan--Reid
\cite[Theorem~9.5.5]{MacReid}).
Concretely, $D$ can also be defined by the Hilbert symbol $(-1,3)$, in
addition to $(-1,-3)$; the quaternion algebra defined over $\Q$ with
symbol $(-1,3)$ gives the totally geodesic surface.  The presence of
an immersed totally geodesic surface implies not just
\fullref{vpbn}, but the stronger statement that $M_0$ has a
finite cover $N$ where $\pi_1(N)$ surjects onto a free group of rank 2
(see Lubotzky \cite{LubotzkyTrans}).  One can also see such a virtually free
quotient directly from the topological description of $M_1$ given in
\fullref{fig:link1}.  There, the underlying space of the orbifold
$M_1$ is $\RP^3 \, \# \,\, \RP^3 \, \# \, L(4,1)$, and hence
$\pi_1(M_1)$ surjects onto $\Z/2\Z * \Z/2\Z * \Z/4\Z$.  The latter group acts
on an infinite $4$--valent tree without a global fixed point.  For an
$M_n$ with $n > 1$, the covering map $M_1 \gets M_n$ gives an action of
$\pi_1(M_n)$ on the same tree, again without a global fixed point.
Therefore, if $M_n$ is a manifold, there is an incompressible surface
dual to this action.  Thus $M_n$ is Haken, and hence all the examples
in \fullref{thm:main-top} are Haken.

\section{Modular forms for $\GL_2/K$}\label{sec:modular-forms}

Let $K/\Q$ be a number field.  We will use $\A_K$ to denote the ad\`eles
of $K$, and $\Af_K$ the finite ad\`eles.  According to the Langlands
philosophy, any regular algebraic cuspidal automorphic representation
$\pi$ of $\GL_2(\A_K)$ should be attached to a compatible family of
$\Gal(\Kbar/K)$--representations, whose local representations are
well behaved and can \emph{a priori} be determined by the local
factors of $\pi$.  If $K$ is totally real such representations may be
constructed ``geometrically'' from the Tate modules of certain Shimura
varieties (see Carayol \cite{Carayol} and Taylor \cite{taylorhilbert}).  We are interested in the case
where $K$ is an imaginary quadratic field, and here the corresponding
symmetric spaces fail to be algebraic varieties.  This causes a great
headache in the construction of Galois representations which one
expects are always geometric.  This problem was solved in the
paper~\cite{taylorothers} by Harris, Soudry and Taylor and the
subsequent paper of Taylor~\cite{taylorimag}. The kernel of the idea
is to use automorphic induction to $\GSp_4/\Q$ where one has geometry
with which to construct Galois representations.  The original desired
$2$--dimensional $\Gal(\Kbar/K)$--representation should then be related
to these $\Gal(\Qbar/\Q)$ representations by Galois induction.  The
result of these labours is the following theorem~\cite{taylorimag}.

\begin{theorem}[Taylor] Let $K$ be an imaginary
  quadratic field, and let $\overline{\phantom{z}}$ denote its non-trivial
  automorphism.  Let $\pi$ be a cuspidal automorphic representation of
  $\GL_2(\A_K)$ such that $\pi_{\infty}$ has Langlands parameter
\[
W_{\C} = \C^{\times} \to \GL_2(\C) \quad \mbox{given by} \quad z \mapsto \left(\begin{matrix} z^{1-k} & 0 \\ 0 & \overline{z}^{1-k}
 \end{matrix}\right), \quad  \mbox{where $k \in \Z_{\geq 2}$.}
\]
Let $S$ be a set of places containing those where $K$ ramifies and
those where $\pi$ or $\pibar$ is ramified.  For $v \notin S$ let
$\{\alpha_v,\beta_v\}$ be the Langlands parameters of $\pi_v$. Let $F$ be the
field generated by $\alpha_v + \beta_v$ and $\alpha_v \beta_v$ for $v \notin S$; it is a
number field. Assume
\begin{enumerate}
\item the central character $\chi$ of $\pi$ satisfies $\chi = \osl{\chi}$,
\item the integer $k$ in the Langlands parameter is even.
\end{enumerate}
Then there is an extension $E/F$ such that for each prime $\lambda$ of $F$
there is a continuous irreducible representation
$$\rho_{\lambda}\co \Gal(\Kbar/K) \to \GL_2(E_{\lambda})$$
such that if $v$ is a
prime of $K$ outside $S$ and not dividing the residue characteristic
$\ell$ of $\lambda$ then $\rho$ is unramified. Moreover, the characteristic
polynomial of $\rho_{\lambda}(\Frob_{v})$ is $(x - \alpha_v)(x - \beta_v)$ for a set
of places $v$ of density one.
\label{theorem:taylor}
\end{theorem}

The existence of Galois representations in a setting where the
associated symmetric spaces are not algebraic now allows us to study
these spaces, which are arithmetic hyperbolic manifolds, using Galois
representations.  Nevertheless, \fullref{theorem:taylor} is not
sufficient for our purposes.  In order to study more precisely the
arithmetic of $\rho$, we need finer control over the behavior of $\rho$ at
primes in the set $S$, and the primes dividing the residue
characteristic $\ell$. What we need can be described by the following
conjectural extension of \fullref{theorem:taylor}:

\begin{conj}[Langlands for $\GL_2/K$] 
  Let $K$ be an imaginary quadratic field, and let $\pi$ be as in
  \fullref{theorem:taylor}, without assuming conditions $1$ or
  $2$.  Then the representation $\rho_{\lambda}$ exists as above, and is
  potentially semistable at $v$ for all $v$.  Furthermore, the
  associated representation of the Weil--Deligne group satisfies the
  local Langlands correspondence with the associated representation
  $\pi_{v}$. 
\label{conj:global}
\end{conj}

One may ask how far Taylor's theorem is from establishing
\fullref{conj:global}. The condition on the central character
seems unavoidable in the arguments of Taylor et al.  Beyond this,
there are essentially two issues, one minor, one more serious.  The
first concerns the behavior at a prime $v \nmid \ell$ when $\pi_v$ and
$\pibar_v$ have different ramification behavior, for example when
$\pibar_v$ is ramified and $\pi_v$ is not. In this context, Taylor's
result does not guarantee (as \fullref{conj:global} does) that
$\rho_{\lambda}$ is unramified at $v$.  (Note the statement of
\cite[Theorem~A]{taylorothers} contains a somewhat mischievous misprint:
``does not divide $\mathfrak{n} \ell$'' should read ``does not divide
$N_{K/\Q}(\mathfrak{n}) \ell$''.)  The issue is that the associated
form $\twpi$ on $\GSp_4$ couples $\pi_v$ and $\pibar_v$ together, and
thus $\twpi_{\ell}$ can be ramified even if $\pi_v$ is not.  One possible
approach would be to exploit the compatibility of local and global
Langlands for $\GSp_4/\Q$; Taylor has sketched an argument for the
first author along these lines. However, such an argument would
require a proof of said compatibility for $\GSp_4/\Q$, which is not
currently known (though not expected to be \emph{too} difficult, for
example it is now known for $\GL_n$; see Taylor--Yoshida \cite{TT}).  A more serious issue
arises for $v | \ell$, even for $k=2$, the case of interest. Here
it is not even known that the $\ell$--adic representations on $\GSp_4$
are potentially semistable at $\ell$, since their construction uses
congruences to higher weight forms. An analogous issue arises for
Hilbert modular forms and was solved also by
Taylor~\cite{taylorhilbert}, but that argument does not generalize to
this case. Nonetheless, we feel that \fullref{conj:global} is a
most reasonable conjecture to make.

The automorphic forms arising in \fullref{conj:global} are
 naturally associated to cohomology classes on their associated
symmetric spaces, which in this case are 3--manifolds. We now
recall the adelic construction of these manifolds. Fix a
quadratic imaginary field $K$, which for ease of
exposition we
assume to have class number one.  Let $\Ok:=\Ok_K$ be the ring of
integers of $K$. Let $\mathfrak{p}$, $\mathfrak{q}$ be two distinct
prime ideals of $\Ok$. Let $U$ be the open compact subgroup of
$\GL_2(\Af_K)$ such that
\begin{enumerate}
\item $U_v$ is $\GL_2(\Ok_v)$ for all $v$
outside $\mathfrak{p}$ and $\mathfrak{q}$.
\item $U_{\mathfrak{p}}$ is the set of matrices in $\GL_2(\Ok_{\mathfrak{p}})$ of the
form
$$\left(\begin{matrix} * & * \\ 0 & * \end{matrix}\right)
\mod \mathfrak{p}.$$
\item $U_{\mathfrak{q}}$ is the set of matrices in $\GL_2(\Ok_{\mathfrak{q}})$ of the
form
$$\left(\begin{matrix} 1 & * \\ 0 & 1 \end{matrix}\right)
\mod {\mathfrak{q}}^n.$$
\end{enumerate}
Define
\[
X\left(\Gamma_0(\mathfrak{p}) \cap \Gamma_1(\mathfrak{q}^n)\right) = 
\raisebox{-0.7ex}{$\GL_2(K)$}\big \backslash \raisebox{0.5ex}{$\left( \GL_2(\Af_K) /U \right) \times \H$}
\]
Now let $D/K$ be the quaternion algebra ramified exactly at
$\{\mathfrak{p},\mathfrak{q}\}$. Let $B$ be a maximal order of $D$.
Let $\m$ be a maximal bi-ideal of $B$ that is trivial away from
$\mathfrak{q}$.  Let $V$ be the compact subgroup of $\GL_2(\Af_K)$ that
is $(1 + \m^n)_v$ for each $v$ (note, this will equal $B^{\times}_v$
for all $v \neq \mathfrak{q}$).  Define
\[
X[\mathfrak{p} \mathfrak{q}^n] =
\raisebox{-0.7ex}{$D^{\times}$} \big \backslash \raisebox{0.5ex}{$\left(
    \GL_2(\Af_K) \big / V \right) \times \H$}
\]
The 3--orbifolds
$X\left(\Gamma_0(\mathfrak{p}) \cap \Gamma_1(\mathfrak{q}^n)\right)$ and
$X[\mathfrak{p} \mathfrak{q}^n]$ have more concrete descriptions as
arithmetic quotients of $\H$; the former as $\H/(U \cap \GL_2(\Ok_K))$
and the latter as 
\[
\H/\left(B^{\times} \cap (1 + \m^n)\right).
\]
There is a well known correspondence between the cohomology of the
3--manifolds constructed above and automorphic forms (see for
example Taylor \cite[Section~4]{taylorimag} or Harder~\cite{Harder}).  A classical version
of this association relates $H^1_\cusp$ of modular curves to classical
modular cusp forms.  Recall here that for a manifold $X$ with
boundary, the cuspidal cohomology $H^1_\cusp(X;\C)$ is the kernel of
the natural homomorphism
$$H^1(X;\C) \to H^1(\partial X;\C).$$  
The cuspidal cohomology of $X\left(\Gamma_0(\mathfrak{p}) \cap
\Gamma_1(\mathfrak{q}^n)\right)$ is exactly associated to spaces of automorphic
forms $\pi$ that satisfy the conditions of
\fullref{conj:global}, and thus we may study the cohomology by
studying the associated (predicted) Galois representations.  The
cohomology of the compact manifold $X[\mathfrak{p} \mathfrak{q}^n]$
can also be associated to automorphic forms. The main theorem of
Jacquet--Langlands~\cite[Section~16]{JL} implies that the resulting space
forms will correspond in a precise way to a subset of the automorphic forms
arising from $H^1_\cusp\left(X\left(\Gamma_0(\mathfrak{p}) \cap
\Gamma_1(\mathfrak{q}^n)\right);\C\right)$. Thus we obtain the following result, which
remains completely mysterious from a topological point of view:

\begin{theorem}[Jacquet--Langlands] 
  If $H^1(X[\mathfrak{p}\mathfrak{q}^n];\C)$ is non-zero then
\[
H^1_\cusp\left(X\left(\Gamma_0(\mathfrak{p}) \cap  \Gamma_1(\mathfrak{q}^m)\right);\C\right) \quad \mbox{is non-zero for
  some $m \leq 2 n$.}
\]
\end{theorem}

A more precise version of this theorem relates not only these
cohomology groups as $\C$--vector spaces but also as modules for the
action of the so called Hecke operators (namely, the correspondence
preserves eigenspaces and eigenvalues).  However, we only apply this
result to conclude that the former cohomology group vanishes because
the latter one does, and so we suppress this extra detail. 
In particular, the conclusion of this theorem that we will need is the
following:

\begin{theorem} 
\label{theorem:before}
Assume \fullref{conj:global}.
Suppose that $H^1(X[\mathfrak{p} \mathfrak{q}^n];\Q)
\neq 0$. Suppose that $\mathfrak{p}$ has residual
characteristic $p$. Then there exists
a field $[E:\Q_p] < \infty$ and a continuous
Galois representation
$\rho\co \Gal(\Kbar/K) \to \GL_2(E)$
such that
\begin{enumerate}
\item $\rho$ is ordinary at $\mathfrak{p}$,
\item $\rho$ is unramified outside $\mathfrak{q}$ and primes dividing
$p$.
\item $\det(\rho) = \chi \psi$, where $\psi$ is a
finite character unramified outside $\mathfrak{q}$.
\item $\rho$ is irreducible.
\end{enumerate}
\end{theorem}

\begin{proof} By the universal coefficient theorem, if $\Q$--cohomology is
  trivial then so is $\C$--cohomology; hence $H^1(X[\mathfrak{p}
  \mathfrak{q}^n];\C) \neq 0$.  By Jacquet--Langlands, the non-triviality
  of this cohomology implies that
  $H^1_\cusp\left(X\left(\Gamma_0(\mathfrak{p}) \cap
      \Gamma_1(\mathfrak{q}^m)\right);\C\right)$ is also 
      nonzero, for some $m$.  As
  mentioned above, the non-vanishing of this cuspidal cohomology gives
  an automorphic representations $\pi$ which satisfies the hypotheses of
  \fullref{conj:global}.  The Galois representation produced by
  this conjecture then satisfies (1)--(4).
\end{proof}

\fullref{thm:main-arith} may be restated in the language of this section
as follows:

\begin{theorem} Assume \fullref{conj:global} and the GRH.  
  Let $K = \Q(\sqrt{-2})$, and let $3 = \pi \pibar$ in $\Ok_K$.  Then
  $H^1(X[\pibar \pi^n]; \Q) = 0$ for all $n$.  
\label{theorem:main}
\end{theorem}

Taking $\mathfrak{p} = \pibar$ and $\mathfrak{q} = \pi$ in
\fullref{theorem:before}, we see that if $X[\pibar \pi^n]$ has
nontrivial $H^1$ then there exists a $3$--adic Galois representations
satisfying properties (1)--(4).  Thus it suffices to prove that no
such representations exist.  This is exactly the content of
\fullref{theorem:main1}, where it is shown that such a Galois
representation satisfying (1)--(3) must violate (4).  We complete the
proof of this purely arithmetic result in the next two sections of
this paper.

\section{Residual representations of small Serre level}\label{sec:Serre-level}

As a warm up to proving \fullref{theorem:main1}, we first study
the possible residual parts of the Galois representations in question.
Let $\F$ be a finite field of residue characteristic $p$. Let $K$ be a
number field, and
$$\rhobar\co \Gal(\Kbar/K) \to \GL_2(\F)$$
a semisimple Galois
representation of Serre conductor $N(\rhobar)$.  A theorem of
Tate~\cite{Tate} shows that if $(K,p,N(\rhobar)) = (\Q,2,1)$ then
$\rhobar$ is trivial.  We prove the following variant of Tate's
results over an imaginary quadratic field.

\begin{theorem}
\label{theorem:3}
 Suppose that $(K,p,N(\rhobar)) =
(\Q(\sqrt{-2}),3,1)$.  Let $3 = \pi \pibar$ and suppose
furthermore that
 $\det(\rhobar) = \chi$, the cyclotomic character restricted to $K$.
Then assuming the GRH,    $\rhobar = \chi \oplus 1$.
\end{theorem}

\begin{proof} Without loss of generality we  assume that $\rhobar$ cannot
be conjugated  so that its image lies in $\GL_2(\F')$ for some proper subfield
$\F' \subset \F$.
We shall break the proof up into
 three cases. Either
\begin{enumerate}
\item The representation $\rhobar$ is reducible.
\item The representation $\rhobar$ is solvable but not reducible.
\item The image of $\rhobar$ is non-solvable.
\end{enumerate}

If $\rhobar$ is reducible and  semisimple it breaks up
as the direct sum of two characters:
$$\eta, \eta'\co \Gal(\Kbar/K) \to \F^{\times}.$$
If $L$ is the maximal abelian
extension of $K$ unramified outside $3$ of order coprime to 3, then by
class field theory $\Gal(L/K)$ is isomorphic to $(\Ok/3 \Ok)^{\times} \simeq
\F^{\times}_3 \times \F^{\times}_3$ modulo $\Ok^{\times} = \{\pm 1\}$. Since this quotient
has order $2$, it follows that $L = K(\zeta_3)$.  Since $\det \rhobar =
\chi$ it follows that $\eta$ and $\eta'$ are (possibly after reordering) the
trivial and cyclotomic character respectively.

If $\rhobar$ is solvable and irreducible then the image of
$$\rhobarproj\co \Gal(\Kbar/K) \to \mathrm{PGL}_2(\F)$$
is either
dihedral, $A_4$ or $S_4$.  Since $\det \rhobar = \chi$ the image of
$\det \rhobarproj$ is $\F^{\times}_3/\F^{\times 2}_3 \simeq \Z/2\Z$.  Since $A_4$
has no such quotients it follows that the image of $\rhobarproj$ is
either dihedral or $S_4$.  Let $M/K$ be the field cut out by the
kernel of $\rhobarproj$.

We first consider the case where $M/K$ is dihedral. Suppose the degree
$[M:K]$ is divisible by $3$. Then the image of $\rhobar$ contains an
element of order $3$, and it follows that the image of $\rhobar$ lands
inside a Borel subgroup of $\GL_2(\F)$ (and is thus reducible) or
contains $\SL_2(\F')$ for some subfield $\F' \subset \F$. The former case
has already been considered, and in the latter case the projective
image $\rhobarproj$ is either non-solvable or $S_4$.  Thus $3$ does
not divide the order of $[M:K]$.  Let $L$ be the maximal extension of
$K$ inside $M$ such that $L/K$ is abelian.  Since $[L:K]$ is coprime
to $3$, we find (as in the reducible case) that $\Gal(L/K)$ is a
quotient of the cokernel $\{ \pm 1\} \to (\Ok/3 \Ok)^{\times} \simeq \F^{\times}_3 \times
\F^{\times}_3$.  Since this cokernel has order $2$, it follows that $L \subseteq
K(\zeta_3) = \Q\left(\sqrt{-2},\sqrt{-3}\right)$. The extension $M/L$ is
abelian since $\rhobarproj$ has dihedral image.  Suppose that $M/L$
was ramified at some prime above $3$. Then $M/K$ is totally ramified
at this prime, and thus the inertia group of $\Gal(M/K)$ is
non-abelian, and has order coprime to $3$.  Yet the maximal tame
quotient of inertia is pro-cyclic, and thus $M/L$ must be unramified
everywhere.  Since $\CL\left(\Q\left(\sqrt{-2},\sqrt{-3}\right)\right)
= 1$, this is impossible.

Now assume that $M/K$ is an $S_4$--extension. Let $L$ be the maximal
extension of $K$ inside $M$ such that $L/K$ is abelian. Since
$\Gal(M/K) = S_4$, it follows that $\Gal(L/K) = \Z/2\Z$, and
$([L:K],3) = 1$. Thus as above $L = K(\zeta_3) =
\Q\left(\sqrt{-2},\sqrt{-3}\right)$.  Let $J$ be the maximal abelian
extension of $L$ contained in $M$. Then $\Gal(J/K) = S_3$. The only
$S_3$ extensions of $K = \Q(\sqrt{-2})$ unramified outside $3$ and
containing $L$ are $L\left(\sqrt[3]{\pi}\right)$,
$L\big(\sqrt[3]{\pibar}\big)$ and $L\left(\sqrt[3]{3}\right)$. All
these extensions are totally ramified over $K$ at both primes above
$3$. Consider the inertia subgroup at a prime above $3$ of
$\Gal(M/K)$. It is a subgroup of $S_4$ that surjects onto $S_3$, and
thus is either $S_3$ or $S_4$.  Since $S_4$ cannot occur as the
inertia group of an extension of $\Q_3$, it follows that $M/L$ is
unramified at all primes above $3$ and thus unramified everywhere. Yet
the class number of each of the three prospective fields $L$ is one,
and thus $M$ does not exist.

Finally, let us assume that the image of $\rhobar$ is non-solvable.
We will need the following lemma.
\begin{lemma}
\label{lemma:Serre}
 Let $\F$ be a finite field
of characteristic $p$, and
let
$$\rho\co \Gal(\Kbar/K) \to \GL_2(\F)$$
be a continuous Galois representation. Fix
a prime $\mathfrak{p}$ above $p$ in $K$. Then either
\begin{enumerate}
\item The image of inertia at $\mathfrak{p}$ is tame.
\item The image of the decomposition group
at $\mathfrak{p}$ is reducible.
\end{enumerate}
\end{lemma}
The proof of the lemma is a standard application of the fact that
$p$--groups do not act freely on finite dimensional $\F_p$--vector
spaces, and we leave it to the reader.  

Let us consider the restriction of $\rhobar$ to inertia at $\pi$.  If
the image of inertia under $\rhobar$ has order coprime to $3$, then
the largest power of $\pi$ dividing $\Delta_{L/K}$ is $\pi^{[L:K] - 1}$. The
contribution to the root discriminant $\delta_{L/\Q} =
|\Delta_{L/\Q}|^{1/[L:\Q]}$ is thus bounded above by $3^{1/2}$. Suppose
the image of the decomposition group at $\pi$ is reducible. All
characters $\Gal(\Qbar_3/\Q_3) \to \F^{\times}$ are the product of an
unramified character and a power of the cyclotomic character. Thus the
tame quotient has order $2$, and the image of wild inertia is
contained in the set of matrices of the form
$$\left(\begin{matrix} 1 & * \\ 0 & 1 \end{matrix}\right) \subset
\GL_2(\F).$$
This subgroup is isomorphic to $\F$. In particular, it is
an elementary three group of order $3^m = \ \|\F\|$.  An explicit
calculation (see Tate~\cite{Tate}) now shows that this is the largest
possible power of $\pi$ dividing $\Delta_{L/K}$ is $(2 + 1/3 - 1/2 \cdot
3^{1-m})[L:K]$.  A similar calculation also applies to $\pibar$.  Thus
if $3^m = \|\F\|$, then the root discriminant $\delta_{L/\Q}$ is bounded
as follows:
\begin{equation}\label{eq:delta}
\log(\delta_{L/\Q}) \leq \log(\delta_{K/\Q}) + \log(3) \left(2 + 1/3 - 1/2 \cdot 3^{1-m}\right).
\end{equation}
On the other hand, since $\det \rhobar = \chi$, we know that $L$
contains $\Q\left(\sqrt{-2},\sqrt{-3}\right)$.  Let $G =
\Gal(L/K)$ and $H =
\Gal\left(L/\Q\left(\sqrt{-2},\sqrt{-3}\right)\right)$.  Since $G$ is
non-solvable the image of $G$ in $\PGL_2(\F)$ is either $A_5$ or
contains the group $\PSL_2(\F)$ (see Serre \cite[Lemma~2]{Serre}).  Since
$\det \rhobar = \chi$ the image of $\det \rhobarproj$ is $\F^{\times}_3/\F^{\times
  2}_3 \simeq \Z/2\Z$.  Since $A_5$ admits no such quotient we see that
the image of $G$ must be all of $\PGL_2(\F)$. It follows that the
image of $H$ is all of $\PSL_2(\F)$, and thus by a classification of
subgroups of $\SL_2(\F)$ that $H = \SL_2(\F)$.  In particular,
\begin{equation}\label{eq:degree}
[L:\Q] = \#H \cdot [\Q\left(\sqrt{-2},\sqrt{-3}\right):\Q]  =  4 \| \SL_2(\F) \| = 4 (3^{2m} - 1) 3^{m}.
\end{equation}
On the other hand, the GRH discriminant bounds of
Odlyzko~\cite{odlyzko} imply that for sufficiently small $\delta_{L/\Q}$
the degree $[L:\Q]$ is bounded above.  Letting $B$ denote this upper
bound for various values of $\delta_{L/\Q}$ we have the following table:

\begin{center}{%
\renewcommand{\arraystretch}{1.4}
\begin{tabular}{c|c|c|c}
$m$ & Upper bound on $\delta_{L/\Q}$ from \eqref{eq:delta} & $B$ &   $4 \| \SL_2(\F) \| $ \\
\hline
$2$ &  $2^{3/2} 3^{13/6} = 30.5708639321$ &  $2400$  & $2880$\\
$3$ &  $2^{3/2} 3^{41/18} = 34.5399086640$ & $10000$  & $78624$\\
$4 \ldots \infty$  & $2^{3/2} 3^{7/3}
=36.7136802477$ & $100000$ & $\geq 2125440$\\
\end{tabular}
}\end{center} For each possible value of $m$, we have $[L : \Q] = 4 \|
\SL_2(\F) \| > B \geq [L : \Q]$, a contradiction.  It follows (on the
GRH) that $\rhobar$ must have solvable image, and thus by previous
considerations must be $\chi \oplus 1$.
\end{proof}

\begin{remark}
  We note that there is a plausible technique for removing the
  assumption on the GRH. If we assume that $\rhobar$ is \emph{modular}
  (which is sufficient for our applications), then a generalized
  Serre's conjecture would imply that $\rhobar$ arises from a mod--$3$
  modular form of sufficiently small level and weight which can then
  be computed explicitly. However, we prefer to assume the GRH
  conjecture, since although it is probably more difficult, it has the
  benefit of already being widely considered.  In contrast, possible
  generalizations of level lowering that may as yet reveal unexpected
  phenomena.
\end{remark}

\section{Global representations of small conductor}\label{sec:global}

As discussed at the end of \fullref{sec:modular-forms}, in order
to complete the proof of \fullref{thm:main-arith} it suffices to
show:

\begin{theorem}
 Let $K = \Q(\sqrt{-2})$,
and let $3 =\pi \pibar$ in $\Ok_K$.
Let $E$ be a local field of mixed residue characteristic
$3$. Let 
$$\rho\co \Gal(\Kbar/K) \to \GL_2(E)$$
be a continuous Galois representation such that
\begin{enumerate}
\item $\rho$ is ordinary at $\pibar$.
\item $\rho$ is unramified outside $\pi$ and $\pibar$.
\item $\det(\rho) = \chi \psi$, where $\psi$ is a
finite character unramified outside $\pi$.
\end{enumerate}
Then $\rho$ is reducible.
\label{theorem:main1}
\end{theorem}

\begin{proof}
  By choosing a lattice in $E$ we may assume that $\rho$ has an integral
  representation $V$.  Note that such a choice is not necessarily
  unique.  Let $\F = \Ok_E/\mathfrak{m}_E$ be the (finite) residue
  field of characteristic three. The fixed field of $\osl{\psi}$ is an
  extension of $K$ of degree coprime to $3$ unramified outside $\pi$.
  By class field theory, such extensions are classified by
  $(\Ok/\pi \Ok)^{\times} \simeq \F^{\times}_3$ modulo global
  units.  Since the image of $-1$ generates $\F^{\times}_3$, it follows
  that $\psi$ is trivial modulo $3$.  Thus $\rhobar$ has determinant
  $\chi$, and hence by \fullref{theorem:3}, the semisimplification
  $\rhobar^{ss}$ is isomorphic to $\chi \oplus 1$.  Thus $V$ has a
  filtration by modules isomorphic to $\Z/3\Z$ or $\mu_3$. By
  assumption $\rho$ is ordinary at $\pibar$. Thus there exists a
  filtration
  $$0 \to V' \to V \to V'' \to 0$$
  of $\Gal(\Kbar_{\pi}/K_{\pi}) =
  \Gal(\Qbar_3/\Q_3)$--modules where $V' \otimes \Q$ and $V'' \otimes \Q$ are
  free of rank one over $E$, and the filtered pieces of $V'/p^n$
  (respectively $V''/p^n$) are all isomorphic as
  $\Gal(\Qbar_3/\Q_3)$--modules to $\mu_3$ (resp.~$\Z/3\Z$).  It
  therefore suffices to show that $V'/p^n$ and $V''/p^n$ extend to
  $\Gal(\Kbar/K)$--modules.  Assume otherwise. Then there must exist an
  extension class in $\Ext^1(\mu_3,\Z/3\Z)$ that splits completely at
  $\pibar$.  It thus suffices to prove the following:

\begin{lemma} 
  The group of extensions $\Ext^1(\mu_3,\Z/3\Z)$ that are unramified
  outside $\pi$ and split completely at $\pibar$ are trivial.
\label{lemma:schoof3}
\end{lemma}

\begin{proof}
  Any such Galois extension defines a $\chi^{-1} = \chi$ extension of
  $K(\zeta_3)$.  Such extensions are of the form $L = K(\zeta_3,\gamma^{1/3})$
  where $\gamma \in K$.  Moreover, $L/K$ is totally split at $\pibar$ if
  and only if one can take $\gamma \equiv 1 \mod \pibar^2$.  Since $L/K$ is
  unramified outside $\pi$ and $\CL(\Ok_K) = 1$ it must be the case
  that $\gamma = \pm \pi$. Yet $\pm \pi \equiv \pm 4 \mod \pibar^2$, and thus this
  extension does not split completely at $\pibar$ (indeed it is
  ramified at $\pibar$).\end{proof}

Having established the lemma, we've proven the theorem, and, as a
consequence, completed the proof of
\fullref{thm:main-arith}.
\end{proof}

We note the following corollary of Theorems~\ref{theorem:main1} 
and~\ref{theorem:3}.
\begin{cor} Assuming the GRH, there does not exist
an abelian variety of $\GL_2$--type over $K = \Q(\sqrt{-2})$ with good
reduction outside the prime $\pi = 1 -  \sqrt{-2}$.
\end{cor}

\begin{proof} 
  Consider such a variety $A$, and let $\rho$ be an associated
  $\GL_2$--representation occurring inside the $3$--adic Tate module of
  $A$. By \fullref{theorem:3} the representation $\rhobar$ is
  reducible. Since $\Z/3\Z$ and $\mu_3$ are the only finite flat group
  schemes of order $3$ over $\Spec(\Ok[1/\pi])$ it follows that $A$ has
  good ordinary reduction at $\pibar$.  From
  \fullref{theorem:main1} we conclude that $\rho$ is reducible.  It
  follows that $A$ has CM by some order in $\Ok$. In particular, $A$
  arises from base change from some totally real field, and thus from
  $\Q$. This implies that $A$ has good reduction at $\pi$, and thus
  good reduction everywhere over $\Ok$.  The result then follows
  from Schoof \cite{schoof}.
\end{proof}

\section{The Boston--Ellenberg approach: pro--$\!p$ groups}\label{sec-pro-p}

As mentioned in the introduction, Boston and Ellenberg
\cite{BostonEllenberg} were able to improve
\fullref{thm:main-arith} to an unconditional statement, not
dependent on the GRH and \fullref{conj:global}; their approach
was to analyze the $M_n$ using the theory of pro--$\!p$ groups rather
than automorphic forms and Galois representations.  In this section,
we will generalize their approach slightly so that it applies, as we
illustrate, to a range of examples, both arithmetic and non-arithmetic.

In order to state the main result, we first need to give a number of
definitions.  To start, if $M$ is a hyperbolic 3--manifold, we say a
tower of finite covers
\[
  M_0 \gets M_1 \gets M_2 \gets M_3 \gets \cdots 
\]
\emph{exhausts} $M$ if $\injrad(M_n) \to \infty$ as $n \to \infty$.  If $\beta_1(M_n)
= H^1(M_n; \Q) = 0$ for all $n$, then we say $M$ can be
\emph{exhausted by rational homology spheres}.  On a different note,
any hyperbolic 3--manifold $M$, arithmetic or not, has an associated
quaternion algebra $A_0(M)$, see eg Maclachlan--Reid \cite[Chapter~3.3]{MacReid}.  Here
we are using the non-invariant quaternion algebra, which can change a
little under finite covers.
 
Turning now to finite groups, the crucial definition is:

\begin{definition}
  Let $p$ be an odd prime.  A finite $p$--group $S$ is \emph{powerful}
  if $S/S^p$ is abelian, where $S^p$ is the subgroup generated by all
  $p^{\mathrm{th}}$ powers.  For $p = 2$, $S$ is powerful if $S/S^4$
  is abelian.
\end{definition}

One should think of powerful $p$--groups as close to being abelian, and
hence sharing many of the properties of abelian groups.  For general
groups, the following definition is central to our criterion:

\begin{definition}
  Let $G$ be a finitely generated group.  We say that $G$ is
  \emph{$p$--powerful} if every $p$--group quotient of $G$ is powerful.
\end{definition}

As we will explain in \fullref{rmk-check-pow}, it is
straightforward to check from a presentation whether $G$ is
$p$--powerful.  The main result of this section is:

\begin{theorem}\label{thm-3-mfld-pro-p}
  Let $M$ be a hyperbolic 3--manifold whose quaternion algebra ramifies
  at a prime of norm $p^n$ where $p \in \Z$ is prime.  Suppose that
  $\beta_1(M) = 0$ and that $|H_1(M;\Z)|$ is coprime
  to $p^{2n} - 1$.  If $\pi_1(M)$ is $p$--powerful then $M$ can be
  exhausted by rational homology spheres.
\end{theorem}

We will actually prove somewhat more than this, but the above is the
one that is easy to apply in practice (see
\fullref{sec-pro-p-examples}).  It follows easily from the
next proposition; we thank Alex Lubotzky for pointing out this
formulation, which removes some unnecessary restrictions in our
original version.

\begin{proposition}\label{prop-pro-p}
  Let $G$ be a finitely generated group which is $p$--powerful.  If
  $\beta_1(G) = 0$, then any $H \unlhd G$ of $p$--power index also has $\beta_1(H) =
  0$.
\end{proposition}

The proof of the proposition uses some fairly elementary aspects of
the theory of pro--$\!p$ groups.  However, in some cases one can reduce
\fullref{thm-3-mfld-pro-p} down to a single non-trivial theorem
about finite $p$--groups; we give that argument in
\fullref{subsec-alt-approach}.  We refer the reader to the book
by Dixon, du Sautoy, Mann and Segal \cite{DixonEtAl1999} for the needed
background about pro--$\!p$ groups.

\begin{proof}[Proof of \fullref{prop-pro-p}]
  For this proof, it is natural to replace $G$ and $H$ by their
  pro--$\!p$ completions $\widehat{G}$ and $\widehat{H}$; since $[G : H]
  = p^n$, the completion $\widehat{H}$ is a finite index open subgroup
  of $\widehat{G}$.  Note $G$ has a surjective homomorphism  onto $\Z$ if
  and only if $\widehat{G}$ has one to $\Z_p$, and the same for $H$.
  Thus it is enough to show that if $\widehat{H} \twoheadrightarrow
  \Z_p$ then so does $\widehat{G}$.
  
  For the rest of this proof, we work exclusively with pro--$\!p$ groups,
  and so drop the decorations from $\widehat{G}$ and $\widehat{H}$ and
  denote them simply $G$ and $H$.   Now suppose that $H
  \twoheadrightarrow \Z_p$.  It suffices to show that some quotient of
  $G$ surjects onto $\Z_p$, and we first exploit this to reduce to the
  case where $G$ is \emph{uniformly} powerful and $H$ is isomorphic to
  $\Z_p^d$, for $d > 0$.
  
  Consider the maximal torsion-free abelian quotient of $H$:
  \[
  1 \to  K \to  H \to  \Z_p^d \to  0.
  \]
  As $K$ is characteristic in $H$, it is normal in $G$, and so we can replace
  $(G, H)$ with $(G/K, H/K = \Z_p^d)$.  As $G$ is powerful, the set $T$ of
  all torsion elements in $G$ is in fact a finite subgroup and $G/T$ is
  uniformly powerful; thus replacing $(G, H)$ with $(G/T, H/(T \cap H) \cong H)$
  reduces to the desired case.
  
  Now $H = \Z_p^d$ is a uniformly powerful open subgroup of the
  uniformly powerful $G$.  Uniformly powerful groups are in fact
  $p$--adic analytic groups, and we will make use of their associated
  Lie algebras.  (For another pro--$\!p$ approach requiring less
  machinery, see \cite[Chapter~4, Example~9]{DixonEtAl1999}.)  Since $H$ is
  abelian, the Lie bracket on $L(H)$ is trivial.  The inclusion
  induced monomorphism of the Lie algebras $L(H) \to L(G)$ then implies
  that the Lie bracket on $L(G)$ is trivial as well.  This forces $G$
  to be abelian and hence isomorphic to $\Z_p^d$
  \cite[Corollary~7.16]{DixonEtAl1999}.  Thus $G$ surjects onto $\Z_p$, as
  desired.
\end{proof}

\begin{proof}[Proof of \fullref{thm-3-mfld-pro-p}]
  Let $G = \pi_1(M)$.  By \fullref{prop-pro-p}, every normal
  subgroup of $G$ of $p$--power index has $\beta_1 = 0$.  Thus to complete
  the proof, we need to show that $G$ is residually a $p$--group, or
  equivalently, that $M$ can be exhausted by regular covers of degree
  $p^n$.  If $\mathfrak p$ is the given prime where $A_0(M)$ ramifies,
  we can consider the principal congruence covers $M({\mathfrak p}^n)$
  of $M$ of level $\mathfrak p^n$ just as in \fullref{subsec:M_n}.
  The first of these covers $M \gets M({\mathfrak p})$ is cyclic of
  degree dividing $p^{2n} - 1$ (see Maclachlan--Reid \cite[Theorem~6.4.3]{MacReid} for
  further details).  The constraint on $H_1(M; \Z)$ forces this
  cover to be trivial.  The remaining covers $M({\mathfrak p}^n) \gets
  M({\mathfrak p}^{n+1})$ all have degree dividing $p^2$.  Just as
  before, these are regular covers of $M$ which exhaust it, completing
  the proof.
\end{proof}

\begin{remark}\label{rmk-check-pow}
  To apply \fullref{thm-3-mfld-pro-p}, we need to be able to check
  that a given finitely presented group is $p$--powerful.  First, for a
  finite $p$--group $S$, consider the lower exponent--$p$ central series
  \[
  S = P_0(S) \rhd P_1(S) \rhd \cdots \rhd P_k(S) = \left\langle  1 \right\rangle  \quad \mbox{where $P_{i+1}(S) = P_i(S)^p [ P_i(S), S]$.}
  \]
  Here the successive quotients are isomorphic to $(\Z/p\Z)^n$ and $k$
  is called the exponent--$p$ class of $S$.  Any finitely presented
  group $G$ has a maximal $p$--quotient of exponent--$p$ class $k$; see
  Newman--O'Brien \cite{NewmanOBrien1996} for details and GAP
  \cite{GAP4} or Magma \cite{Magma}
  for an implementation.  Now suppose $G$ is a finitely presented
  group, and $S$ the maximal $p$--quotient of exponent--$p$ class $2$.
  We claim that $G$ is $p$--powerful if and only if $S$ is powerful;
  this follows by noting that if $Q$ is a non-powerful $p$--quotient of
  $G$, then $Q/P_2(Q)$ is a non-powerful quotient of $S$ (see the
  lemma below).  Since $S$ is computable from a presentation from $G$
  and checking if $S$ is powerful is easy since it is finite, we have
  the needed test.

  We turn now to the $p$--group lemma just used. 
  \begin{lemma}
    Let $Q$ be a $p$--group.  Then $Q$ is powerful if and only if
    $Q/P_2(Q)$ is powerful.  
  \end{lemma}
  \begin{proof}  
    First suppose $p$ is odd.  Noting that $Q/P_2(Q)$ is powerful if
    and only if
    \[
    (Q/Q^p)/P_2(Q/Q^p)
    \]
    is abelian, it suffices to show
    that if $T$ has exponent $p$ then $T$ is abelian if and only if
    $T/P_2(T)$ is.  If $T/P_2(T)$ is abelian then $P_2(T) \geq [T, T] =
    P_1(T)$ as $T^p = 1$.  Thus the series stabilizes at $P_1(T) =
    [T,T] = 1$, and $T$ is abelian, as needed.
    
    If $p = 2$, we can similarly reduce to showing that a group $T$ with
    exponent $4$ is abelian if $T/P_2(T)$ is abelian.  We proceed by
    induction on $\abs{T}$.  If $T$ has exponent $2$, it is abelian so
    assume that $T^2$ is a non-trivial subgroup.  Then as $T$ is
    2--group, $N = T^2 \cap Z(T)$ is non-trivial.  By induction, $G/N$ is
    abelian, and hence $N \geq [T, T]$.  Thus commutators are all central
    which implies $[t_1^2, t_2^2] = [t_1, t_2]^4 = 1$ for $t_i \in T$;
    thus $T^2$ is abelian of exponent 2.  Then $[T, T] \leq T^2$ is as well,
    which means $[t_1^2, t_2] = [t_1, t_2]^2 = 1$.  Thus $P_2(T) = 1$
    and $T$ is abelian as desired.
  \end{proof}
\end{remark}

\subsection{Examples}\label{sec-pro-p-examples}

While the hypotheses of \fullref{thm-3-mfld-pro-p} are easy to
check, one does not expect to find so many examples.  To say that
$\pi_1(M)$ is $p$--powerful implies that pro--$\!p$ completion is analytic,
and in fact these two conditions are equivalent, up to passing to
subgroups of finite index.  Lubotzky has shown that if $\dim H_1(M;
\F_p) > 3$ then the pro--$\!p$ completion of $\pi_1(M)$ is \emph{not}
analytic (see Lubotzky \cite{Lubotzky1983}, and apply Theorem~1.2 as
strengthened by Remark~1.4, to the combination of Lemma~5.1 and
Lemma~1.1).  In fact, because a pro--$\!p$ group generated by two
elements is solvable and $\pi_1(M)$ is not, it is easy to see that
\fullref{thm-3-mfld-pro-p} can only apply when $\dim H_1(M; \F_p)
= 3$.

One example where \fullref{thm-3-mfld-pro-p} applies is the
manifold $M_2$ from \fullref{subsec:M_n}, and this was Boston and
Ellenberg's original insight \cite{BostonEllenberg}.  We searched the
Hodgson--Weeks census \cite{SnapPea} for manifolds where
\fullref{thm-3-mfld-pro-p} applies. First, in order for $\pi_1(M)$
to have a  non-cyclic pro--$\!p$ completion, it is necessary that $\dim
H_1(M; \F_p) \geq 2$.  Among the 11{,}126 manifolds in the census there
are 719 such pairs $(M, p)$ where $p$ is an odd prime, typically $3$
or $5$.  (Technical aside: we exclude $v1539(5,1)$ which has $H_1(M;
\Z) = \Z^2$ and hence has this property for all $p$.)  Then we looked
at the lower exponent--$p$ central series of $G = \pi_1(M)$ as discussed
in \fullref{sec-pro-p-examples}.  If we set $d_k = \dim
P_k(G)/P_{k-1}(G)$ we observed an apparent dichotomy of behavior:
\begin{enumerate}
  \item The $d_k$ are uniformly bounded by $3$. (61\% of the cases.)
  \item The $d_k$ eventually start to grow exponentially with $k$. (39\% of the cases.)  
\end{enumerate}
\begin{figure}[ht!]
\labellist\small
\pinlabel {$2$} [t] at 51 0
\pinlabel {$4$} [t] at 94 0
\pinlabel {$6$} [t] at 137 0
\pinlabel {$8$} [t] at 181 0
\pinlabel {$10$} [t] at 224 0
\pinlabel {$k$} [l] at 275 0
\pinlabel {$2$} [r] at 0 12
\pinlabel {$3$} [r] at 0 51
\pinlabel {$4$} [r] at 0 89
\pinlabel {$5$} [r] at 0 127
\pinlabel {$6$} [r] at 0 165
\pinlabel {$\log(d_k)$} [b] at 0 170
\endlabellist
\centerline{\includegraphics[scale=0.85]{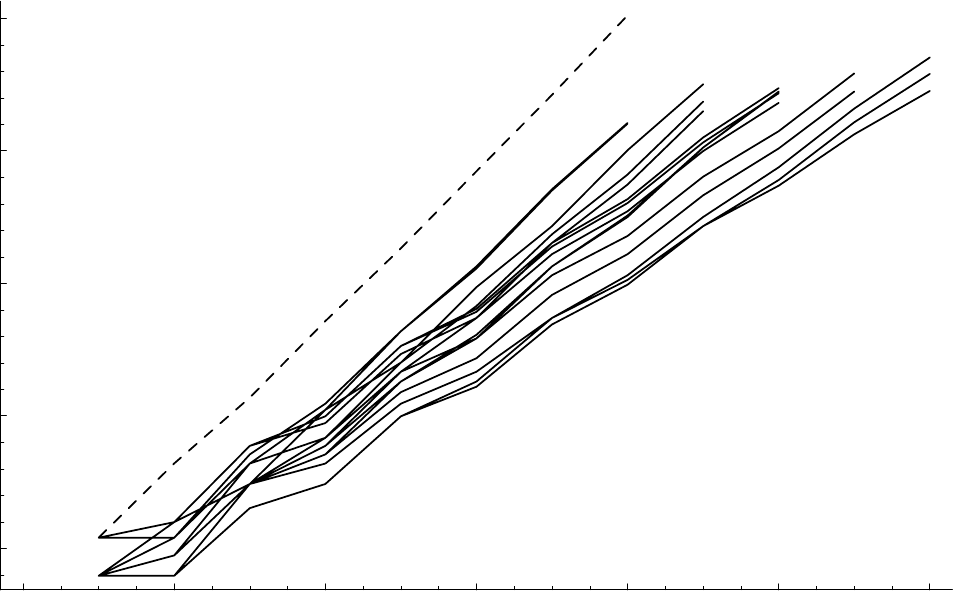}}\vspace{10pt}
\caption{This figure shows the exponential growth of the $d_k$.  The solid
lines plot $\left(k - k_0, \log (d_{k-k_0}) \right)$ for each group,
where $k_0$ is the first $k$ for which $d_k > 5$ (that is, the first
few small $d_k$ are omitted to save space).  The dotted line is the same
quantity for the free group of rank 2.  Even though this plot shows 278
different groups, the $d_k$ fall into only 15 distinct patterns.}
\label{fig-exp-growth}
\end{figure}
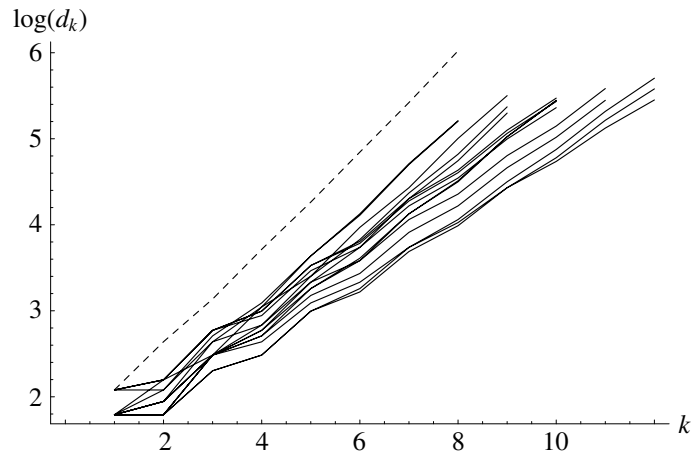
These should correspond precisely to whether the pro--$\!p$ completion of
$G$ is analytic or not. In the former case, we confirmed this by
checking that some $P_k(G)$ is in fact powerful (where $k \leq 2$) in
all but one example when $p \leq 7$, but did not check the 37 examples
where $p > 7$.  In the second case, that these groups must be
non-analytic follows from the result of Lubotzky mentioned above.  The
apparent exponential growth is clear from \fullref{fig-exp-growth}.
Most of these groups have 2 generators, and so for comparison it is
interesting to consider the free group of rank 2.  In that case, for
any prime $p$, Witt's formula gives $d_k = \sum_{m = 1}^k (1/m) \sum_{d|m}
\mu(m/d) 2^d$, which implies exponential growth at rate $\log 2$; that
is $\lim_{k \to \infty} (1/k) \log d_k = \log 2$.  As you can see from
\fullref{fig-exp-growth}, the non-analytic groups seem quite close
to the free group in this sense.

Turning now to those manifolds where the pro--$\!p$ completion is
analytic, we need to check the quaternion algebra hypothesis of
\fullref{thm-3-mfld-pro-p}, which can be done using the program
Snap \cite{SnapArticle,Snap}.  Unfortunately, Snap often fails to
compute the needed information for the larger examples, so we content
ourselves with simply listing 20 examples where
\fullref{thm-3-mfld-pro-p} applies, as opposed to an exhaustive
list.  In all cases, the prime at which the quaternion algebra ramifies
has norm a rational prime.
\begin{equation*}
\begin{split}
p &= 3 : 
\underline{s649(-5,1)}, {m007(3,2)}, \underline{s784(3,1)}, \underline{s784(-5,1)}, {v2715(-1,2)}, \\
& \quad \quad \quad{v3215(1,2)}, \underline{v3215(5,1)} \\
p &= 5: {m003(-3,1)}, {m285(1,2)}, \underline{v2616(4,1)}, \underline{v2674(-1,3)}, \underline{v2933(1,4)}, \\
&\quad \quad \quad \underline{m216(4,3)}, {s889(3,2)}, {v2739(1,2)} \\ 
p &= 7: {m017(-3,2)}, \underline{m017(-5,1)}, \underline{m022(-5,1)}, \underline{m030(1,3)}, {s527(-3,2)} 
\end{split}
\end{equation*}
Here the non-arithmetic examples are underlined.  With the exception of
the last two examples for $p=3$, the theorem is really applied to the
cover corresponding to the kernel of $\pi_1(M) \to H_1(M ; \F_p)$.  Also
note that the first example for $p=5$ is the Weeks manifold.  

We also looked at what happened for $p = 2$.  There are $1492$ census
manifolds where $\dim H^1(M ; \F_2) \geq 2$.  The situation was similar
to before, but a new behavior appears:
\begin{enumerate}
\item The maximal $2$--quotient of $\pi_1(M)$ is finite.  (34\% of the cases.)
\item The pro--$2$ completion of $\pi_1(M)$ is analytic and apparently
  infinite. (28\% of the cases.)
\item The pro--$2$ completion of $\pi_1(M)$ is non-analytic and the
  $d_i$ appear to have an exponential growth rate $(0.25, 0.55)$.
  (39\% of the cases.)
\end{enumerate}
Here are 19 examples where \fullref{thm-3-mfld-pro-p} applies to the kernel
of $\pi_1(M) \to H_1(M ; \F_2)$: 

\begin{flushleft}
$m039(6,1)$, $m035(-6,1)$, $m037(2,3)$, $s227(-2,3)$, $s961(1,2)$,
$s781(-5,1)$, $s786(4,1)$, $v2229(-4,3)$, $v2231(5,1)$, $v2230(4,3)$,
$s957(-1,4)$, $s955(-1,4)$, $s961(1,4)$, $v3273(5,1)$, $s594(3,2)$,
$s956(4,1)$, $s961(4,1)$, $s960(-1,4)$, $v3111(2,3)$.
\end{flushleft}

Here the first three are arithmetic, and all the others are non-arithmetic.
Complete software (in Magma \cite{Magma}) and data files for all of
the above may be obtained at \cite{paperwebsite}.

\subsection{Alternate approach}\label{subsec-alt-approach}

It is possible to circumvent most of the pro--$\!p$ machinery used in the
proof of \fullref{prop-pro-p}, at least if one increases the
hypothesis somewhat.  In this section we give a version which is
enough to apply to many of the examples in the preceding section, yet
relies on only a single theorem about $p$--groups.  

Throughout we fix a prime $p$.  For any group, let $d(S) = \dim(H^1(S;
\F_p))$; for a $p$--group, this is equal to the minimal number of
generators (see Dixon--du Sautoy--Mann--Segal
\cite[Theorem~0.9]{DixonEtAl1999}).  Then one has the following
basic fact about powerful $p$--groups, which is one of the ways they
are analogous to abelian groups.

\begin{theorem}[Dixon--du Sautoy--Mann--Segal {{\cite[page
41]{DixonEtAl1999}}}] Let $S$ be a powerful $p$--group.  If $H$ is any subgroup of $S$ then
  $d(H) \leq d(S)$.
\end{theorem}

Using this we will show

\begin{proposition}\label{prop-pro-p-alt}
  Let $G$ be a finitely presented group.  Suppose $G$ is $p$--powerful
  for $p \geq 5$, and that $d(G) \leq 3$ and $\beta_1(G) = 0$.  If $N$ is
  a normal subgroup of $G$ of index $p^n$ then $\beta_1(N)$ = 0.
\end{proposition}

As noted in \fullref{sec-pro-p-examples}, in when
\fullref{thm-3-mfld-pro-p} applies one always has $d(\pi_1(M)) \leq
3$.  Thus \fullref{prop-pro-p-alt} can stand in for
\fullref{prop-pro-p} in the proof of
\fullref{thm-3-mfld-pro-p} whenever $p \geq 5$.

\begin{proof}  First, we claim that $\beta_1(N) \leq  3$; this follows if we can show
  $d(N) \leq 3$.  Consider the mod--$p$ abelianization map
  \[
  1 \to  K \to  N \to  H_1(N; \F_p) \to  0.
  \]
  Then as $K$ is a characteristic subgroup of $N$, it is normal in $G$.  By
  assumption, the $p$--group $G/K$ is powerful.  Thus by the above theorem
  \[
  d(N) = d(N/K) \leq  d(G/K) \leq  d(G) \leq 3
  \]
  as desired.
 
  To complete the proof, it suffices to show that $\beta_1(N) = \beta_1(G)$.
  Consider the action of $S = G/N$ on $H_1(N; \Z)/(\mathrm{torsion})$.
  This gives us a homomorphism
  \[
  S \to  \GL_{ \beta_1(N)} \Z  \leq  \GL_3 \Z.
  \]
  If this map is trivial, then one has an isomorphism between $H_1(N;
  \Q)$ and $H_1(G; \Q)$ as desired.  But it is easy to see this map is
  trivial --- just note that the only primes which are orders of elements
  of $\GL_3 \Z$ are 2 and 3, since the characteristic polynomial of an
  element of $\GL_3\Z$ has degree 3.
\end{proof}

\section{Twist-knot orbifolds}\label{sec:twist}

In this section, we investigate the congruence covers of twist-knot
orbifolds.  In contrast with the rest of the paper, our approach here
is experimental --- we simply examine large numbers of such covers.
Interestingly, we find an apparent dichotomy of behavior between the
arithmetic and non-arithmetic examples.  In particular, the congruence
covers of the arithmetic orbifolds were much more likely to have $\beta_1
> 0$.

\begin{figure}[ht!]
\labellist\small
\pinlabel {$\frac1n$} [l] at 265 450
\endlabellist
\centerline{\includegraphics[scale=0.6]{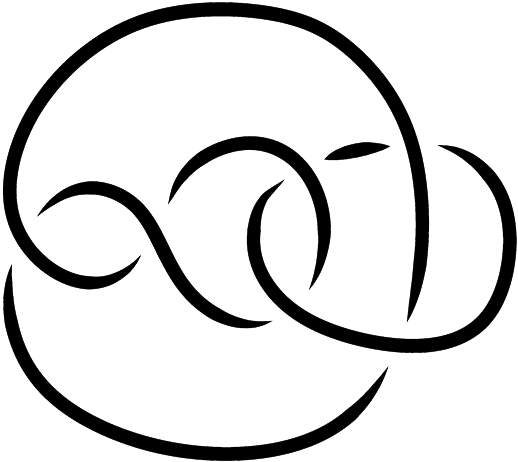}}
\caption{$K_n$ is
  obtained by doing $\frac1n$ Dehn surgery on the
  Whitehead link.}
\label{fig:twist}
\end{figure}

Let us first roughly outline the results; a precise account follows.
The twist knots $K_n$, for $n \in \Z$, are a simple family of knots
described in \fullref{fig:twist}.  Here, $K_1$ is the trefoil,
$K_0$ the unknot, and $K_{-1}$ the figure--8.  We define $T(n, k)$ to
be the orbifold whose underlying space is $S^3$, and singular set
consists of the twist knot $K_n$ labeled by $k$.  We focused on the
$T(n,k)$ which are hyperbolic and where $\abs{n} \leq 4$ and $3 \leq k \leq
7$; we examined $31$ of these orbifolds.  For each $T(n,k)$, we looked
at congruence covers of the form $\Gamma_0(\mathfrak{p})$, which are
covers of degree $N(\mathfrak{p}) + 1$.  We did this for all
$\mathfrak{p}$ of norm less than $10{,}000$, which gave about
600--1200 covers for each $T(n,k)$.  The natural question is: what
percentage of these covers have $\beta_1 > 0$?  Keeping the base orbifold
$T(n,k)$ fixed, we observed two distinct kinds of behavior
\begin{enumerate}
  \item A large proportion, 26--66\% of the covers $\Gamma_0(\mathfrak{p})$ have $\beta_1 > 0$. 
  \item Less than 2\% of the covers have $\beta_1 > 0$.
\end{enumerate}
About a third of the $T(n,k)$ fell into the first category, with the
rest having the rare $\beta_1 > 0$ behavior.  From a topological point of
view, this dichotomy is rather odd since these manifolds are closely
related (see in particular the formula \eqref{eq:twist-pi1} for
$\pi_1$).  Surprisingly, this dichotomy corresponds to exactly whether
$T(n,k)$ is arithmetic, with the arithmetic manifolds falling into
case (1).  This table summarizes the data:

\begin{table}[h]
\definecolor{agray}{gray}{0.8}
\begin{center}
\begin{tabular}{rr|*{7}{c}}
\multicolumn{2}{c}{} & \multicolumn{7}{c}{Twist parameter $n$ in $K_n$} \\
 & & $-4$ & $-3$ & $-2$ & $-1$ & $2$ & $3$ & $4$ \\
\cline{2-9}\\[-12.9pt]
& 3 &  0.9 & \cellcolor{agray} 49.0 & \cellcolor{agray} 55.5 &  & \cellcolor{agray} 40.8 & \cellcolor{agray} 41.6 & 1.3   \\
& 4 &   1.0 & 0.8 & \cellcolor{agray} 55.0 & \cellcolor{agray} 40.7 & \cellcolor{agray} 65.5 & 1.5 & 0.5   \\
\raisebox{1.5ex}[0pt]{Orbifold} & 5   & 0.8 & 0.7 & 0.8 & \cellcolor{agray} 54.8 & \cellcolor{agray} 56.8 & 0.9 & 0.9   \\
\raisebox{1.5ex}[0pt]{label $k$} & 6  & 1.1 & 0.8 & 0.7 & \cellcolor{agray} 36.3 & \cellcolor{agray} 26.0 & 0.7 & 1.6   \\
& 7 &   0.9  & 0.9 & 0.7 & 2.0 & 1.7 & 0.4 & 1.1   \\
\end{tabular}
\end{center}
\caption{Percentage of $\Gamma_0(\mathfrak{p})$ with $\beta_1 > 0$, where $N(\mathfrak{p}) \leq 10{,}000.$}
\end{table}
The shaded boxes are the arithmetic $T(n,k)$.  Here, $T(0,k)$ and
$T(-1,3)$ were omitted because they are not hyperbolic.  The number of
covers examined was in the range $[578,632]$ when $k$ was even, and in
$[1168,1247]$ when $k$ was odd.  Which $T(n,k)$ are arithmetic was
determined by Hilden--Lozano--Montesinos-Amilibia \cite{HLM}.  The arithmetic $T(n,k)$ all fall into the
classes of such known to satisfy \fullref{vpbn}.  In a few
cases, this is because the trace field has a subfield of index 2 and
so the results of Labesse--Schwermer \cite{LabesseSchwermer} and Lubotzky
\cite{LubotzkyArith} applied, but in most cases
one needs to use the result of Clozel \cite{Clozel}.

For the non-arithmetic examples, it is worth listing the norms of the
exceptional primes $\mathfrak{p}$ where $\beta_1\left( \Gamma_0(\mathfrak{p})
\right) > 0$:
 
{\centering\small
\begin{tabular}{ll}
$T(-4, 3)$ & $157, 197, 239, 3^5, 2^8, 257, 271, 293, 349, 7507$\\
$T(4, 3)$ & $13, 29, 41, 53, 61, 73, 89, 5^3, 127, 2^7, 151, 173, 233, 587, 1201$\\
$T(-4, 4)$ & $73, 79, 103, 233, 1999$ \\
$T(-3, 4)$ & $103, 113, 11^2, 167, 7759$\\
$T(3, 4)$ & $31, 41, 79, 97, 103, 137, 167$\\
$T(4, 4)$ & $23, 103$ \\
$T(-4, 5)$ & $2^6, 79, 101, 131, 13^2, 191, 239, 241$ \\ 
$T(-3, 5)$ & $19, 29, 7^2, 71, 79, 139, 311, 761$ \\ 
$T(-2, 5)$ & $29, 7^2, 71, 79, 89, 2^8, 379, 5^4$\\ 
$T(3, 5)$ & $7^2, 59, 61, 71, 79, 101, 131, 13^2, 271, 881$ \\ 
$T(4, 5)$ & $11, 19, 59, 61, 89, 11^2, 131, 239, 1361, 2099, 4049$ \\
$T(-4, 6)$ & $5^2, 59, 97, 107, 181, 2^8, 6659$ \\
$T(-3, 6)$ & $11, 61, 71, 349, 2053$ \\
$T(-2, 6)$ & $5^2, 157, 181, 937$ \\
$T(3, 6)$ & $23, 71, 647, 5209$ \\
$T(4, 6)$ & $23, 47, 71, 83, 97, 2^7, 131, 229$ \\ 
$T(-4, 7)$ & $13, 3^3, 29, 41, 43, 97, 127, 449, 1093, 2633$\\
T$(-3, 7)$ & $13, 29, 41, 43, 7^2, 113, 139, 1483$\\
$T(-2, 7)$ & $13, 3^3, 127, 181, 503, 7^4$ \\ 
$T(-1, 7)$ & $2^6, 83, 113, 13^2, 181, 211, 239, 3^6, 29^2, 41^2, 43^2, 71^2, 97^2$ \\ 
$T(2, 7)$ & $13, 29, 41, 7^2, 83, 97, 113, 127, 139, 181, 349, 463, 6007$ \\ 
$T(3, 7)$ & $13, 113, 211, 307, 617$ \\
$T(4, 7)$ & $29, 41, 43, 97, 127, 139, 379, 1721$
\end{tabular}
}

We examined two of the $T(n,k)$ more closely, computing additional
covers out to $N(\mathfrak{p}) \leq 25{,}000$.  The first of these is
$T(4,4)$, which above has just 2 covers with $\beta_1 > 0$.  This pattern
continued; there were 727 covers $\Gamma_0(\mathfrak{p})$ with $N(\mathfrak{p})$ in
$[10{,}000, 25{,}000]$, and none of these had $\beta_1 > 0$.  Should it
be possible to answer Cooper's question from a topological or
combinatorial viewpoint, $T(4,4)$ might be a good place to start.
Indeed, it seems plausible that $T(4,4)$ has only finitely many
congruence covers of the form $\Gamma_0(\mathfrak{p})$ which have $\beta_1 >
0$.  Perhaps this is true of the principal congruence covers as well.
In which case, trying to prove \fullref{vpbn} solely by looking
an congruence covers would be a bad idea.  The other one we looked at is
$T(-1,7)$.  In this case, the trend again continues, with $\beta_1 \left(
  \Gamma_0(\mathfrak{p}) \right) > 0$ for a series of primes
$\mathfrak{p}$ with square norm: $113^2, 127^2, 139^2$.  Complete
data, including the Magma source code used for the computations, can
be found at the accompanying website for this paper \cite{paperwebsite}.

\subsection{Details}

We now describe the $T(n,k)$ in more detail, explain exactly which
covers we looked at, and give a few computational caveats.  The
orbifold fundamental group of $T(n,k)$ is given by:
\begin{equation}\label{eq:twist-pi1}
\Gamma = \bigl\langle a,b ~\big|~ a^k = b^k = 1, w^n a = b w^n \bigr\rangle
\quad\text{where $w = b a^{-1} b^{-1} a$.}
\end{equation}
Assuming $T(n,k)$ is hyperbolic, consider the representation $\rho_0 \maps
\Gamma  \to \PSL_2\C$ corresponding to the hyperbolic structure.
Regardless of whether $T(n,k)$ is arithmetic, Weil's local rigidity
theorem \cite{Weil} shows that $\rho$ can be conjugated so that the
image lies in $\PSL_2 L$, for some number field $L$.    Explicitly, 
\[
\rho_0(a) = \begin{pmatrix} e^{i\pi/k} & 1 \\ 0 & e^{-i\pi/k} \end{pmatrix} \quad %
\rho_0(b) = \begin{pmatrix} e^{i\pi/k} & 0 \\ t & e^{-i\pi/k} \end{pmatrix}
\]
where $t$ is an algebraic integer satisfying a monic polynomial $r(z)$
in $\Z[e^{i \pi/k}]$ of degree $\approx 2 \abs{n}$ (see
Hoste--Shanahan \cite[Section~2]{HosteShanahan05}).  

Let $L = \Q(e^{i \pi/k}, t)$ be the field generated by the entries of
$\setdef{\rho_0(\gamma)}{\gamma \in \Gamma}$, and let $K = \Q( e^{i \pi/k} + e^{-i
  \pi/k}, t)$ be the field generated by $\tr(\rho_0(\gamma))$.  Then $[L :
K]$ is at most $2$, and indeed is 2 in the cases at hand.  Let
$\mathfrak{p}$ be a prime ideal of $\Ok_K$, and let $\mathfrak{q}$ be
a prime ideal of $\Ok_L$ which divides $\mathfrak{p}$.  Now consider
the congruence quotient $\Gamma \to \PSL_2(\Ok_L/\mathfrak{q})$.  By
definition, the congruence cover $\Gamma_0(\mathfrak{p})$ is the inverse
image of a Borel subgroup under this map --- that is, those $\gamma \in \Gamma$
which are upper-triangular mod $\mathfrak{q}$.  The reason
$\Gamma_0(\mathfrak{p})$ can be thought of as a function of $\mathfrak{p}$
rather than $\mathfrak{q}$ is that if the choice for $\mathfrak{q}$ is
not unique, ie $\mathfrak{p}$ splits in $\Ok_L$, then the two
choices for $\mathfrak{q}$ give congruence quotients with the same
kernel.  For all but finitely many $\mathfrak{p}$ which do split
in $\Ok_L$, the homomorphism $\Gamma \to \PSL_2(\Ok_L/\mathfrak{q})$ is
surjective; in this case, $\Gamma_0(\mathfrak{p})$ has index $\|
P^1(\Ok_L/\mathfrak{q}) \| = N(\mathfrak{p}) + 1$.  For all but
finitely many $\mathfrak{p}$ which do not split in $\Ok_L$, the
homomorphism $\Gamma \to \PSL_2(\Ok_L/\mathfrak{q})$ maps onto a conjugate
of the smaller group $\PSL_2(\Ok_K/\mathfrak{p})$.  Again,
$\Gamma_0(\mathfrak{p})$ has index $N(\mathfrak{p}) + 1$.  (A more elegant
point of view here would be construct these quotients by localizing
the quaternion algebra associated to $\rho(\Gamma)$; this is also the
connection between the definition here and the one described in
\fullref{subsec:M_n}.)

The reason we did not look at the \emph{principal} congruence cover,
which is the kernel of $\Gamma \to \PSL_2(\Ok_L/\mathfrak{q})$, is purely
pragmatic; because those covers are so much larger it's not possible
to get beyond $N(\mathfrak{p})$ a few hundred, and so there would not
be enough data to draw interesting conclusions.  As it was, the
computations took several CPU-months.

We did not implement the above definition of $\Gamma_0(\mathfrak{p})$
directly because of the degrees of some of the fields involved.
Instead, for each finite field $\F_q$, we searched directly for
\emph{epimorphisms} $\Gamma \to \PSL_2\F_q$, in a way that generates all
the congruence quotients, as well as a possibly a few additional
epimorphisms for small $\F_q$.  Following
Hilden--Lozano--Montesinos-Amilibia \cite{HLM} and Hoste--Shanahan
\cite{HosteShanahan05}, an irreducible representation $\pi_1(T(n,k))
\to \PSL_2 E$, for any field $E$, is essentially determined by $\tr(
\rho_0(a)^2)$ and $\tr( \rho_0(ab) )$.  These quantities must satisfy
certain integer polynomials, and conversely any solution gives
a representation.  Thus we simply searched for solutions to these
equations over $\F_q$ to find the needed epimorphisms.  However, one
must be careful as the ``trace variety'' defined by these equations is
not always irreducible over $\Q$.  When it is reducible, we worked out
the subvariety containing the image of the canonical representation
$\rho_0$, and then used the equations defining that subvariety in our
search.  This ensures that only finitely many non-congruence covers are
generated.  When counting epimorphisms, we considered two epimorphisms
$\Gamma \to \PSL_2\F_q$ equivalent if they different by an automorphism
of $\PSL_2(\F_q)$; this differs from counting the number of ideals of
$\Ok_K$ of norm $q$.

Finally, for speed reasons, when we determined whether $\beta_1 > 0$ we
cheated a bit and worked over the finite field $\F_{31991}$ rather than
$\Q$.  Thus there could be some false positives where the cover really
doesn't have $\beta_1 > 0$, though of course there are no false negatives.
Based on the experience of Dunfield--Thurston \cite{DunfieldThurston}, we
expect there are at most a handful of such false positives, if any at all.

\bibliographystyle{gtart}
\bibliography{link}

\end{document}